\documentclass[12pt,twoside]{amsart}
\usepackage{hyperref}
\usepackage{amsthm, amsmath, amscd, amssymb,centernot}
\usepackage[all]{xy}
\usepackage[T1]{fontenc}
\usepackage[left=2cm,top=2.5cm,bottom=3cm,right=3cm]{geometry}

\usepackage{amsthm}
\usepackage{tikz-cd}
\usepackage{amssymb}
\usepackage{enumerate}
\usepackage{amsmath} 
\usepackage[mathscr]{euscript}
\makeatletter
\@namedef{subjectclassname@2020}{%
\textup{2020} Mathematics Subject Classification}
\makeatother
\newtheorem{thm}{Theorem}[section]
\newtheorem{lem}[thm]{Lemma}
\newtheorem{prop}[thm]{Proposition}
\newtheorem{defi}[thm]{Definition}

\newtheorem{corl}[thm]{Corollary}

\theoremstyle{plain}
\numberwithin{equation}{section}

\newtheorem*{theorem*}{Theorem}

\theoremstyle{definition}

\begin{document}

\begin{abstract}
We describe the GIT compactification of the moduli of (2,2)-type effective divisors of $\mathbb{P}^1\times\mathbb{P}^2$ (i.e.,  surfaces of the linear system $\vert \pi_1^*\mathcal{O}_{\mathbb{P}^1}(2)\otimes \pi_2^*\mathcal{O}_{\mathbb{P}^2}(2)\vert$ ) which are generically Del Pezzo surfaces of degree two. In order to get the compactification, we characterize stable and semi-stable (2,2)-type surfaces, and also determine the equivalence classes of strictly semi-stable (2,2)-type surfaces. Moreover, we describe the boundary of the moduli of (2,2)-type surfaces.
\end{abstract}

\title{Stability and semi-stability of  $(2,2)$-type surfaces}

\author[A.J. Parameswaran]{A.J. Parameswaran}
\address{School of Mathematics, Tata Institute of Fundamental Research, Homi Bhabha Road, Mumbai 400005, India.}
\email{param@math.tifr.res.in}

\author[Nabanita Ray]{Nabanita Ray}
\address{Chennai Mathematical Institute, H1 SIPCOT IT Park, Siruseri, Kelambakkam 603103, 
India.}
\email{nabanitaray2910@gmail.com / nabanitar@cmi.ac.in}

\subjclass[2010]{Primary 14J10, 14J17, 14L24; Secondary 14C20, 14D20.}
\thanks{}

\date{September, 2020}
\maketitle
\tableofcontents
\section{Introduction}
Mumford's geometric invariant theory provides a construction of moduli spaces of family of varieties. In this paper, we apply his methods to obtain a description of moduli space of (2,2)-type surfaces of  $\mathbb{P}^1\times\mathbb{P}^2$ and it's GIT compactification. To this extent, we study  stable surfaces, semi-stable surfaces, and  strictly semi-stable surfaces. Many authors studied moduli spaces of hypersurfaces of given degree, and classified stable and semi-stable hypersurfaces in terms of their singularities. Hilbert studied plane curves of degree $\leq$ 6 and cubic surfaces in \cite{H}. Shah provided much more detailed information about sextic plane curves \cite{S1}, and analyzed quartic surfaces \cite{S2}. Recently many results are developed for quintic surfaces,  for instance  one can check \cite{Rana}, \cite{G}, etc. There are few development in higher dimensional hypersurfaces. In particular, the stability  and the moduli of cubic threefolds (resp. fourfolds) are studied in \cite{Y1} and \cite{A}(resp. \cite{Y2} and \cite{L}).  For further references see \cite{Mu}. It can be noticed (from \cite{H}, \cite{S2}, \cite{G}, \cite{Rana}) that smooth surfaces are always stable and stable surfaces have  at most Du Val  A-D-E type isolated singularities depending on their degrees. There are not enough works about the stability of effective divisors in a variety $X$ which is different from the projective space.

In this note, we study (2,2)-type effective divisors of $\mathbb{P}^1\times\mathbb{P}^2$ which are parametrized by $\mathbb{P}(H^0(\pi_1^*\mathcal{O}_{\mathbb{P}^1}(2)\otimes \pi_2^*\mathcal{O}_{\mathbb{P}^2}(2)))\simeq\mathbb{P}^{17}$. Note that these surfaces are not hypersurfaces. It is known from \cite{Ra}, that smooth (2,2)-type surfaces are rational and isomorphic to $\mathbb{P}^2$ blown-up at seven points i.e., Del Pezzo surfaces of degree two. There is a natural linear action on the linear system $\vert \pi_1^*\mathcal{O}_{\mathbb{P}^1}(2)\otimes \pi_2^*\mathcal{O}_{\mathbb{P}^2}(2)\vert$ by the reductive group $G=\text{SL}(2,\mathbb{C})\times \text{SL}(3,\mathbb{C})$. We define the  stability and the semi-stability with respect to this action. Unlike the cases of hypersurfaces, in our situation smooth surfaces  need not be stable  (but always semi-stable, see Lemma \ref{lem2.17})  and also a stable (2,2)-type surfaces can have non-isolated singularities (Proposition \ref{prop3.9}).

In this paper we prove our results  using term by term analysis of the defining equations. As a consequence of  the Hibert-Mumfod criterion, one obtains a one-parameter subgroup $\lambda$  of $G$ that acts linearly on $\mathbb{P}^1\times\mathbb{P}^2$. This enables us to choose coordinates on $\mathbb{P}^1\times\mathbb{P}^2$ and with this coordinate, we obtain the equations for the divisors. This process is repeated for each of the divisors. 

This article is arranged as follows. In Section $2$ we recall the preliminaries about  Hilbert-Mumford criterion for stability and classification of isolated surface singularities. In Section $3$, we describe irreducible  semi-stable surfaces and prove the following result.

\begin{thm}\label{thm1.1}
 Let $S=Z(f)$ be an irreducible (2,2)-type surface. Then $S$ is semi-stable if and only if for all $P=P_1\times P_2\in S_{\text{sing}}$ each of the following conditions is satisfied:
  
  (i) The tangent cone of $S$ at $P$ is not a pull back of a tangent cone of $P_2$ by the map $p_2$.

(ii) If the fibre over $P_1$ is non-reduced, then the generic points of the fibre $p_1^{-1}(P_1)$ are not in the ramification locus of $p_2$.

(iii) The fibre over $P_1$ is reduced but non-irreducible and let $L_1$ and $L_2$ be two components of the fibre over $P_1$. If generic points of $L_2$ are in the ramification locus of $p_2$ and $p_2$ is not finite for some $P=P_1\times P_2\in L_2$ which correspond to sections $\sigma: \mathbb{P}^1\rightarrow \mathbb{P}^1\times P_2\subset S$ of $p_1$, then the corresponding map $\phi_{\sigma}:\mathbb{P}^1\rightarrow \mathbb{P}(T_{\mathbb{P}^2,P_2})$ is  $\phi_{\sigma}(\mathbb{P}^1)\neq p_2(L_2)$  (see Notations and conventions \ref{item7}).
\end{thm}

In Section $4$ the following theorem is proved classifying  irreducible stable surfaces.

\begin{thm}\label{thm1.2}
 Let $S=Z(f)$ be an irreducible semi-stable (2,2)-type surface of $\mathbb{P}^1\times\mathbb{P}^2$. Then $S$ is stable if and only if  $S$ satisfies following properties:
 
 (i) Either $p_2:S\rightarrow \mathbb{P}^2$ is a finite map or  there exist some sections $\sigma:\mathbb{P}^1\rightarrow S$ of $p_1$ which are contracted by $p_2$. Moreover, Im$(\sigma)$ contains at most $A_1$-type singular points and $\phi_{\sigma}:\mathbb{P}^1\rightarrow \mathbb{P}(T_{\mathbb{P}^2,P_2})$ (see \ref{item7}) is non-constant for all such sections.
 
 (ii) Let $P=P_1\times P_2$ be a singular point of $S$ and $P$ is not $A_1$-type. Then the fibre over $P_1$ is reduced.
 \end{thm}
Moreover in Proposition \ref{prop3.9} of Section 4, we prove that some stable surfaces have $A_1$, $A_2$ or $A_3$-type isolated singularities and even have non-isolated singularities.
 
Sections $5$ and $6$ are devoted to classifying irreducible strictly semi-stable surfaces and  semi-stable criterion for 
non-irreducible surfaces respectively (cf. Theorem \ref{thm3.11} and Theorem \ref{thm3.16} respectively). 

Mukai described a set of monomials which generate  non-stable hypersurfaces in $\mathbb{P}^n$ in terms of 
normalized one-parameter subgroups of $SL(n+1)$. We extend this idea in Section $7$. We describe bi-graded monomials of type $(2,2)$ which generate non-stable $(2,2)$-type divisors in $\mathbb{P}^1\times\mathbb{P}^2$.

There are some results related to GIT compactification of Del Pezzo surfaces, e.g.  one can check \cite{M} and \cite{MM} for degree three and degree four Del Pezzo surfaces respectively. In the last Section $8$, we describe the GIT compactification of the moduli of (2,2)-type surfaces (Theorem \ref{thm1.3}) i.e., Del Pezzo surfaces of degree two. We describe the dimension of moduli space and the boundary of the moduli space.
\begin{thm}\label{thm1.3}
The moduli of (2,2)-type surfaces of $\mathbb{P}^1\times\mathbb{P}^2$ has dimension 6. 
The closed subset $\Big (\vert(2,2)\vert^{\text{ss}}//G\Big )\backslash \Big(\vert(2,2)\vert^{\text{s}}/G\Big)$ consists of three rational curves with a common point.
\end{thm}

\subsection{Notations and conventions}
We fix the following notations:
\begin{itemize}
 \item\label{item0} {We denote by $\pi_1:\mathbb{P}^1\times\mathbb{P}^2\rightarrow\mathbb{P}^1$ and $\pi_2:\mathbb{P}^1\times\mathbb{P}^2\rightarrow\mathbb{P}^2$ first and second projection maps respectively. Let $S$ be a surface inside $\mathbb{P}^1\times\mathbb{P}^2$. Then $p_1:S\rightarrow\mathbb{P}^1$ and $p_2:S\rightarrow\mathbb{P}^2$ are restrictions of projection maps respectively.}
 \item\label{item1} {
Let $S=Z(f)$ be a (2,2)-type surface of $\mathbb{P}^1\times\mathbb{P}^2$, passing through the point $P=P_1\times P_2$. Also let that $x_i$'s and $y_i$'s are homogeneous coordinates of $\mathbb{P}^1$ and $\mathbb{P}^2$ respectively. Then $f\in V=\mathbb{C}[x_0,x_1]_2\otimes\mathbb{C}[y_0,y_1,y_2]_2$ the space of homogeneous $(2,2)$-type bi-degree polynomials. In some affine neighbourhood of $P$, $f$ can be written as
$$f(x,y)=\sum_{\alpha=(\alpha_1,\alpha_2)\in\mathbb{N}^2\times\mathbb{N}^3} a_{\alpha}(x-P_1)^{\alpha_1} (y-P_2)^{\alpha_2}\text{;  $| \alpha_1|\leq 2$,  $| \alpha_2|\leq 2$. } $$
This is an expansion of $f$ around $P$, where we use the usual multi-index notation.} 
\item\label{item2}{ We denote the degree $d$ homogeneous part of a polynomial $g(z_1,\cdots, z_n)$ by $g(z_1,\cdots, z_n)_d$.}
\item\label{item4}{ The GIT terminology is that of Mumford et al. \cite{MFK}. For us, unstable means
not semi-stable, non-stable is  for not properly stable, and strictly
semi-stable means semi-stable, but not properly stable.}
\item\label{item5}{ We say that one semi-stable surface degenerates to
another if the second lies in the orbit closure of the first.}
\item\label{item6}{ In this note we write categorical (resp. geometric) quotients in place of good categorical (resp. good geometric) quotients.}
\item\label{item7}{ Let $S$ be a $(2,2)$-type surface and $P=P_1\times P_2\in S$. Assume that $p_2$ is not a finite map at $P_2$. Then there is a section $\sigma:\mathbb{P}^1\rightarrow S$ of $p_1$ such that $\sigma(\mathbb{P}^1)\simeq\mathbb{P}^1\times P_2$ and the section is contracted by the map $p_2$. As $p_1:S\rightarrow \mathbb{P}^1$ is a conic bundle map (see \cite{Ra}), fibres are conic. If such a section $\sigma$ exists, then for all $P_1\in \mathbb{P}^1$, $p_2(p_1^{-1}(P_1))$ is a conic in $\mathbb{P}^2$ which passes through $P_2$. If $\sigma(\mathbb{P}^1)\nsubseteq S_{\text{sing}}$, then generic fibres are smooth at $P_2$. Hence we can define a morphism from an open subset $U\subseteq \mathbb{P}^1$ to the tangent space of $P_2\in\mathbb{P}^2$ corresponding to the section $\sigma$, i.e., $\phi_{\sigma}:U\rightarrow \mathbb{P}(T_{\mathbb{P}^2, P_2})\simeq \mathbb{P}^1$ such that $\phi_{\sigma}(P_1)$ is the tangent line of the conic $p_2(p_1^{-1}(P_1))$ at $P_2$. From the valuation criterion the map $\phi_{\sigma}$ can be extended to $\mathbb{P}^1$. In particular, if $S=Z(f(x,y))$ and $P_2=[1,0,0]$ then $\phi_{\sigma}(P_1)=Z(\frac{\partial g}{\partial y_1}\vert_{(0,0)}y_1+\frac{\partial g}{\partial y_2}\vert_{(0,0)}y_2)$, where $g(y):=f(P_1,y)$.}
\end{itemize}

\section{Preliminaries}
In the whole article for basic algebraic geometric notations and results we follow \cite{Ha} and for geometric invariant theory, we follow \cite{MFK}, \cite{Mu} and \cite{D}.

The main tool for analyzing stability and semi-stability is the \it Hilbert-Mumford Numerical Criterion \rm (see \cite{D} Theorem 9.1, or \cite{M} Theorem 7.3 and Theorem 7.4). 

The one-parameter subgroups (1-PS) $\lambda$ of $G=\text{SL}(2)\times\text{SL}(3)$ which are used to check the Hilbert-Mumford Numerical Criterion are assumed diagonal $t\in \mathbb{C}^*\rightarrow \text{ diag }(t^{r_0},t^{r_1})\times \text{ diag }(t^{s_0},t^{s_1}, t^{s_2})\in \text{SL}(2)\times\text{SL}(3)$
 with the weights satisfying $r_0+r_1=0$ and $s_0+s_1+s_2=0$;  $r_0\leq r_1$ and $s_0\leq s_1\leq s_2$.  We call such a one-parameter subgroups \it normalized\rm. Given a monomial $x_0^{\alpha_0}x_1^{\alpha_1}y_0^{\beta_0}y_1^{\beta_1}y_2^{\beta_2}$ its weight w.r.t. a normalized one-parameter subgroups $\lambda(t)=\text{ diag }(t^{r_0},t^{r_1})\times \text{ diag }(t^{s_0},t^{s_1}, t^{s_2})$ is $\alpha_0r_0+\alpha_1r_1+\beta_0s_0+\beta_1s_1+\beta_2s_2$ i.e., the action of $\lambda$ on monomials are
 $$\Big(\text{ diag }(t^{r_0},t^{r_1})\times \text{ diag }(t^{s_0},t^{s_1}, t^{s_2})\Big)x^{\alpha}y^{\beta}= \Big (\text{ diag }(t^{r_0},t^{r_1})\times \text{ diag }(t^{s_0},t^{s_1}, t^{s_2})\Big ) x_{0}^{\alpha_0}x_1^{\alpha_1}y_0^{\beta_0}y_1^{\beta_1}y_2^{\beta_2}$$
$$=t^{r_0\alpha_0+r_1\alpha_1+s_0\beta_0+s_1\beta_1+s_2\beta_2}x_{0}^{\alpha_0}x_1^{\alpha_1}y_0^{\beta_0}y_1^{\beta_1}y_2^{\beta_2}=t^{r\alpha+s\beta}x^{\alpha}y^{\beta}.$$
 
 Let $S=Z(f)$ be a (2,2)-type surface of $\mathbb{P}^1\times\mathbb{P}^2$ and the normalized one-parameter subgroups $\lambda$  acts on 
 $$f=\sum_{(\alpha,\beta)\in \mathbb{N}^2\times\mathbb{N}^3,\mid\alpha\mid=2,\mid\beta\mid=2}a_{\alpha\beta}x^{\alpha}y^{\beta}.$$
 
The Hilbert-Mumford function (see \cite{MFK}, Definition 2.2]) is
$$\mu(f,\lambda)=\text{ min }\{r\alpha+s\beta\mid a_{\alpha\beta}\neq 0\}.$$

A (2,2)-type polynomial $f$ is stable (resp. semi-stable) if and only if $\mu(f,\lambda)<0$ (resp. $\leq 0$) for allone-parameter subgroups $\lambda$ of  $\text{SL}(2)\times\text{SL}(3)$. We fix coordinates of $\mathbb{P}^1\times\mathbb{P}^2$ and assume that all one-parameter subgroups used are normalized (equivalently, fix a maximal torus $T$ in $\text{SL}(2)\times\text{SL}(3)$ and consider only one-parameter subgroups of $T$).

A one-parameter subgroup of $G$ acts on $V$ by the induced action. Hence the map $G_m\rightarrow V$ is given by $t\rightarrow \lambda(t).f$, for all $f\in V$. If this morphism extends to a morphism $\mathbb{A}^1\rightarrow V$, then the image of the origin is called the limit of $\lambda$ at $f$ as $t\rightarrow 0$, and denoted by $\lim\limits_{t\rightarrow 0} \lambda(t) . f$. So $f$ is not semi-stable (resp. not stable) if and only if there exists a one-parameter subgroups $\lambda$ satisfying $\mu(f,\lambda)>0$ (resp. $\mu(f,\lambda)\geq 0$) or equivalently $\lim\limits_{t\rightarrow 0} \lambda(t). f=0$ (resp. $\lim\limits_{t\rightarrow 0} \lambda(t). f$ exists). 
 
Now we recall some notations and results related to the surface singularity.
 Let $X$ and $Y$ be two varieties over $\mathbb{C}$. Then two points $p\in X$ and $q\in Y$ are \it
analytically isomorphic \rm if there is a $\mathbb{C}$-algebra isomorphism $\widehat{\mathcal{O}}_p\simeq\widehat{\mathcal{O}}_q$.

An affine hypersurface $Z(f)$ has $A_n$-type singularity at $P$ if and only if $P$ is analytically isomorphic to the origin of the affine variety $Z(x^2+y^2+z^{n+1})$. This type of singularity is called Du Val singularity (for more about Du Val singularity see \cite{R}).

 There are some useful techniques to study the type of isolated singularity at a given point. E.g.  Mather and Yau \cite{MY} proved
that two germs of complex analytic hypersurfaces of the same dimension with isolated singularities are biholomorphically equivalent if and only if their moduli algebra
are isomorphic. 

Let $\mathcal{O}_{n+1}$ be  the ring of  germs at the origin of $\mathbb{C}^{n+1}$ and $(U,0)$ be a germ at the origin of a  hypersurface in  $\mathbb{C}^{n+1}$. Let $I(U)$
be the ideal of functions in $\mathcal{O}_{n+1}$  vanishing on $U$, and let $f$ be a generator of
$I(U)$. The ring 
$$A(U)=\mathcal{O}_{n+1}/\Big (f,\frac{\partial f}{\partial z_0},\cdots, \frac{\partial f}{\partial z_n}\Big )$$
is called \it  moduli algebra \rm of $U$.

It is well know that $U\backslash\{0\}$ is non-singular if and only if $A(U)$ is finite dimensional as a $\mathbb{C}$ vector space (see \cite{MY}). Therefore $A(U)$ is infinite dimensional if and only if $U$ has non-isolated singularity at $0$.
\begin{defi}\label{defi2.1}
The germ $(U,0)$ is called quasi-homogeneous if $f(z)$ is quasi-homogeneous i.e., $f\in (\frac{\partial f}{\partial z_0},\cdots, \frac{\partial f}{\partial z_n} )$.
\end{defi}

\begin{thm}\label{thm2.10}
 \rm (\cite{MY}) \it Suppose $(U, 0)$ and $(W, 0)$ are germs of hypersurfaces in $\mathbb{C}^{n+1}$ and $U\backslash \{ 0 \}$
is non-singular. Then the following conditions are equivalent.

(i) $(U, 0)$ is biholomorphically equivalent to $(W, 0)$.

(ii) $A(U)$ is isomorphic to $A(W)$ as a $\mathbb{C}$-algebra.
\end{thm}
Now we recall Corollary 2 of Chapter VII, page-137. \cite{ZS} which will be used in some proofs of this note.
 \begin{thm}\label{thm2.13}
  Let $A$ be an integral domain and let $f^1,f^2, \cdots ,f^m$
be $m$ power series in $A[[X_1, X_2 , \cdots , X_n]]$, $m\leq n$, such that the initial
forms of the $f^i$ are linearly independent linear forms $f^1_1,f^2_1, \cdots ,f^m_1$.
Then the substitution mapping $\phi : g( Y_1, Y_2,\cdots, Y_m)\rightarrow g(f^1,f^2, \cdots ,f^m)$
is an isomorphism of $A[[Y_1, Y_2,\cdots , Y_m ]]$ into $A[[X_1, X_2, \cdots ,X_n]]$.
If, furthermore, $m = n$, $Y_i = X_i$ $i = 1, 2,\cdots , n$, and the determinant of
the coefficients of the linear forms $f^1_1,f^2_1, \cdots ,f^m_1$ is a unit in $A$ (in
particular, if $A$ is a field and the above determinant is $\neq 0$), then $\phi$ is an
automorphism of $A[[X_1, X_2,\cdots , X_n]]$.
\end{thm}
The following result follows easily using the definition of $A_n$-type singularities, Theorem \ref{thm2.10} and Theorem \ref{thm2.13}.
\begin{prop}\label{prop2.14}
  Let $f(z)$ be a holomorphic function
in a neighbourhood of the origin in $\mathbb{C}^3$. Let $(U,0)=(Z(f),0)$. Then $Z(f)$ has

(i) $A_n$-type singularity at the origin if and only if $A(U)\cong \mathbb{C}[[z]]/(z^n)$ for $n\geq1$.

(ii) If $A(U)\cong \mathbb{C}[[z]]$, then $Z(f)$ has non-isolated singularity at $0$.
\end{prop}

\section{Irreducible  semi-stable surfaces}
In this section we study a necessary and sufficient  criterion of irreducible semi-stable surfaces. In Theorem 5.25 \cite{Bu}, Bunnett proved that any bi-degree $(d,e)\in\mathbb{Z}^2$ smooth hypersurface of product of projective space $\mathbb{P}^n\times\mathbb{P}^m$ is semi-stable w.r.t. the action of $\text{SL}(n+1)\times \text{SL}(m+1)$.
In the following lemma, we see that not only smooth (2,2)-type surface of $\mathbb{P}^1\times\mathbb{P}^2$ is semi-stable but  any (2,2)-type surface having at most $A_1$-type singularities is also semi-stable.
\begin{lem}\label{lem2.17}
Let $S$ be a (2,2)-type surface of $\mathbb{P}^1\times\mathbb{P}^2$. If $S=Z(f)$ has at most $A_1$-type singularities then $S$ is a semi-stable surface.
 \begin{proof}
Let $S\in \vert(2,2)\vert$ of $\mathbb{P}^1\times\mathbb{P}^2$ and $S=Z(f)$, $f\in V=\mathbb{C}[x_0,x_1]_2\otimes\mathbb{C}[y_0,y_1,y_2]_2$ the space of homogeneous $(2,2)$-type bi-degree polynomial. Then 
  $$f=\sum_{(\alpha,\beta)\in \mathbb{N}^2\times\mathbb{N}^3,\mid\alpha\mid=2,\mid\beta\mid=2}a_{\alpha\beta}x^{\alpha}y^{\beta}.$$
  We prove that if $S$ is unstable, then $S$ has a singular point which is not $A_1$-type.
 Let $S$ be an irreducible unstable surface. Then there exists a normalized one-parameter subgroup $\lambda: G_m\rightarrow G$ such that $\lim\limits_{t \to 0}\lambda(t).f=0$, where
 $$\lambda(t)=\text{diag}(t^{r_0},t^{r_1})\times \text{diag} (t^{s_0},t^{s_1},t^{s_2}),$$
 $r_0+r_1=0$, $r_0\leq r_1$;   $s_0+s_1+s_2=0$, $s_0\leq s_1\leq s_2$; $r_0+s_0<0$. 
 $$\lambda(t).f=\sum a_{\alpha\beta}(t^{r_0}x_0)^{\alpha_0}(t^{r_1}x_1)^{\alpha_1}(t^{s_0}y_0)^{\beta_0}(t^{s_1}y_1)^{\beta_1}(t^{s_2}y_2)^{\beta_2}=\sum a_{\alpha\beta}t^{r\alpha+s\beta}x^{\alpha}y^{\beta}.$$
 $\lim\limits_{t \to 0} \lambda(t).f=0$, i.e., $r\alpha+s\beta>0$ for all $a_{\alpha\beta}\neq 0$. 

We can write  $f(x_0,x_1,y_0,y_1,y_2)$ explicitly as
 
 $f(x_0,x_1,y_0,y_1,y_2)=x_0^2(a_{11}y_1^2+a_{22}y_2^2+a_{12}y_1y_2+a_{01}y_0y_1+a_{02}y_0y_2+a_{00}y_0^2)+x_0x_1(b_{11}y_1^2+b_{22}y_2^2+b_{01}y_0y_1+b_{02}y_0y_2+b_{12}y_1y_2+b_{00}y_0^2)+x_1^2(c_{00}y_0^2+c_{11}y_1^2+c_{22}y_2^2+c_{01}y_0y_1+c_{02}y_0y_2+c_{12}y_1y_2).$
 
 As $S$ is unstable, the coefficient of $x_0^2y_0^2$, $a_{00}=0$ otherwise $r_0+s_0>0$ which is a contradiction.
 Also, the coefficient of $x_0^2y_0y_1$, $a_{01}=0$ otherwise $2r_0+s_0+s_1>0$ implies $s_0+s_1>0$, therefore $s_2< 0$ which is a contradiction.
Similarly, the coefficient of $x_0x_1y_0^2$, $b_{00}=0$ otherwise $r_0+r_1+2s_0>0$ implies $s_0>0$, which gives a contradiction.
Moreover, the coefficient of $x_0x_1y_0y_1$, $b_{01}=0$ otherwise $r_0+r_1+s_0+s_1> 0$ implies $s_2< 0$ which is a contradiction. 

Now we rewrite the polynomial $f$,
 
 $f(x_0,x_1,y_0,y_1,y_2)=x_0^2(a_{11}y_1^2+a_{22}y_2^2+a_{12}y_1y_2+a_{02}y_0y_2)+x_0x_1(b_{11}y_1^2+b_{22}y_2^2+b_{02}y_0y_2+b_{12}y_1y_2)+x_1^2(c_{00}y_0^2+c_{11}y_1^2+c_{22}y_2^2+c_{01}y_0y_1+c_{02}y_0y_2+c_{12}y_1y_2).$
 
 If the coefficient of $x_0^2y_0y_2$, $a_{02}=0$, then it is clear that $S$ is not smooth at the point $P=[1,0]\times[1,0,0]$. If $a_{02}\neq 0$ then there must exist some monomials  which are not divisible  by $y_2$ i.e.,  coefficients of some  monomials $x_0^{\alpha_0}x_1^{\alpha_1}y_0^{\beta_0}y_1^{\beta_1}$ in $f$ are non-zero. Then at least one of  the coefficient of  monomials $x_0^2y_1^2$, $x_1^2y_0^2$, $x_1^2y_0y_1$, $x_0x_1y_1^2$, and $x_1^2y_1^2$ is non-zero.  As $f$ is unstable, $\alpha_0r_0+\alpha_1r_1+\beta_0s_0+\beta_1s_1>0$.  
  
  Note that $a_{02}\neq 0$ and unstability of $f$ imply the following inequality
  \begin{align}\label{ssequ2.3}
   2r_0+s_0+s_2>0.
  \end{align}
If the coefficient of $x_0^2y_1^2$, $a_{11}\neq 0$ then $2r_0+2s_1>0$. This inequality and the inequality of (\ref{ssequ2.3}) imply $r_0>0$ which is a contradiction.
 
 If the coefficient of $x_1^2y_0^2$, $c_{00}\neq 0$ then $2r_1+2s_0>0$. This inequality and the inequality of (\ref{ssequ2.3}) imply $3s_0+s_2>0$. But note that $s_0+s_1+s_2=0$. Hence  $2s_0-s_1>0$ implies $2s_0>s_1\geq s_0$,   which is a contradiction as $s_0\leq0$.
  
  If the coefficient of $x_1^2y_0y_1$, $c_{01}\neq 0$ then $2r_1+s_0+s_1>0$. This and inequality (\ref{ssequ2.3}) imply $s_0>0$ which is a contradiction.
  
  If the coefficient of $x_0x_1y_1^2$, $b_{11}\neq 0$ then $r_0+r_1+2s_1>0$, i.e., $s_1>0$. This and inequality (\ref{ssequ2.3}) imply $r_0>0$ which is a contradiction.
 
 If the coefficient of $x_1^2y_1^2$, $c_{11}\neq 0$ then $r_1+s_1>0$. This and inequality (\ref{ssequ2.3}) imply $r_0>0$. Hence there is a contradiction.
 
 As $S$ is an irreducible unstable surface, $a_{02}=0$. Then clearly $S$ is not smooth at the point $P=[1,0]\times[1,0,0]$. 
 Then the tangent cone of $S$ at $P$  is $$f(1,x_1,1,y_1,y_2)_2=a_{11}y_1^2+a_{22}y_2^2+a_{12}y_1y_2+b_{02}x_1y_2+c_{00}x_1^2.$$ 
The determinant of the Hessian of $f$ at $P$ is $H=8a_{11}a_{22}c_{00}-2a_{11}b_{02}^2-2a_{12}^2c_{00}$.

If $a_{11}\neq 0$ and $c_{00}\neq 0$, then $2r_0+2s_1>0$ and $2r_1+2s_0>0$ respectively. They imply $s_2<0$ which is a contradiction. So $a_{11}$ and $c_{00}$ are not  simultaneously non-zero.

If $a_{11}=0$, then $f(1,x_1,1,y_1,y_2)_2=a_{22}y_2^2+a_{12}y_1y_2+b_{02}x_1y_2+c_{00}x_1^2$ and $H=-2a_{12}^2c_{00}$.

If $c_{00}=0$, then
$f(1,x_1,1,y_1,y_2)_2=a_{11}y_1^2+a_{22}y_2^2+a_{12}y_1y_2+b_{02}x_1y_2$ and $H=-2a_{11}b_{02}^2$.

If $a_{12}\neq0$ and $c_{00}\neq0$, then we have following inequalities
$2r_0+s_1+s_2>0$, and  $2r_1+2s_0>0$ respectively. Together they imply $r_0>0$, which is a contradiction.
If $a_{11}\neq0$ and $b_{02}\neq0$, then we have following inequalities
$r_0+s_1>0$ and  $r_0+r_1+s_0+s_2>0$. They imply $r_0>0$ which is also a contradiction. 

Eventually we get that the tangent cone of $S$ at $P$  is one of the followings and $H=0$.
\begin{align}\label{equ3.2}
  \begin{cases}
   a_{22}y_2^2+b_{02}x_1y_2+c_{00}x_1^2\\ 
 a_{22}y_2^2+a_{12}y_1y_2+b_{02}x_1y_2\\
 a_{11}y_1^2+a_{22}y_2^2+a_{12}y_1y_2
  \end{cases}
 \end{align}

So if $S$ is unstable, then the determinant of the Hessian at $P$ is zero  i.e., $f$ has other than $A_1$-type singularity at $P$. Hence the result follows.
\end{proof}
 \end{lem}
\begin{prop}\label{prop3.3}
Let $S=Z(f)$ be a $(2,2)$-type surface. If there is a section $\sigma:\mathbb{P}^1\rightarrow S$ of $p_1$  such that the image is contracted by $p_2$, i.e., $p_2(\text{Im}(\sigma))=P_2\in\mathbb{P}^2$ and $\text{Im}(\sigma)\subseteq S_{\text{sing}}$, then $S$ is an unstable surface.
\begin{proof}
 Let us assume that $S$ passes through $P=[1,0]\times [1,0,0]$ and the section of $p_1$ is $\sigma:\mathbb{P}^1\rightarrow \mathbb{P}^1\times [1,0,0]$. Then 
 
 $f(x_0,x_1,y_0,y_1,y_2)=x_0^2(a_{11}y_1^2+a_{22}y_2^2+a_{12}y_1y_2+a_{01}y_0y_1+a_{02}y_0y_2)+x_0x_1(b_{11}y_1^2+b_{22}y_2^2+b_{01}y_0y_1+b_{02}y_0y_2+b_{12}y_1y_2)+x_1^2(c_{11}y_1^2+c_{22}y_2^2+c_{01}y_0y_1+c_{02}y_0y_2+c_{12}y_1y_2).$
 
 From the hypothesis, we have $\text{Im}(\sigma)\subseteq S_{\text{sing}}$. Hence each fibre of $p_1$ is a singular conic. Note that $p_1^{-1}([\alpha_0,\alpha_1])=[\alpha_0,\alpha_1]\times Z(\psi)$ has singularity at $[\alpha_0,\alpha_1]\times [1,0,0]$, where $\psi=y_1^2(\alpha_0^2a_{11}+\alpha_0\alpha_1b_{11}+\alpha_1^2c_{11})+y_2^2(\alpha_0^2a_{22}+\alpha_0\alpha_1b_{22}+\alpha_1^2c_{22})+y_1y_2(\alpha_0^2a_{12}+\alpha_0\alpha_1b_{12}+\alpha_1^2c_{12})$. Therefore 
 
 $f(x_0,x_1,y_0,y_1,y_2)=x_0^2(a_{11}y_1^2+a_{22}y_2^2+a_{12}y_1y_2)+x_0x_1(b_{11}y_1^2+b_{22}y_2^2+b_{12}y_1y_2)+x_1^2(c_{11}y_1^2+c_{22}y_2^2+c_{12}y_1y_2).$
 
 Now $\lambda(t)=\text{diag}(t^{-1},t^{1})\times \text{diag} (t^{-4},t^{2},t^{2})$ satisfies $\lim\limits_{t \to 0}\lambda(t).f=0$. Hence $S$ is unstable.
\end{proof}
\end{prop}
In the following lemma, we see some surfaces having higher order singularities are also semi-stable.
 \begin{lem}\label{lem2.19}
Let $S=Z(f)$ be an irreducible (2,2)-type surface. Then $S$ is a semi-stable surface if the conditions (i), (ii), (iii) in Theorem \ref{thm1.1} hold for all $P=P_1\times P_2\in S_{\text{sing}}$.


\begin{proof}
  It is enough to prove that if 
  $S=Z(f)$ is an irreducible unstable surface, then  some  $P=P_1\times P_2\in S_{\text{sing}}$  satisfy  one of the following properties:
  
  (i) The tangent cone of $S$ at $P$ is a pull back of a tangent cone of $P_2$ by the map $p_2$.
  
  (ii) The fibre over $P_1\in\mathbb{P}^1$ is non-reduced and generic points of $p_1^{-1}(P_1)$ are in the ramification locus of $p_2$.
  
  (iii) The fibre over $P_1$ is reduced and non-irreducible and $L_2$ is   one of the components of the fibre over $P_1$. Then generic points of $L_2$ are in the ramification locus of $p_2$ and $p_2$ is not finite at some $P\in L_2$  which corresponds to  a section $\sigma: \mathbb{P}^1\rightarrow \mathbb{P}^1\times P_2\subset S$ of $p_1$. Moreover, the corresponding map $\phi_{\sigma}:\mathbb{P}^1\rightarrow \mathbb{P}(T_{\mathbb{P}^2,P_2})$ is constant and  $\phi_{\sigma}(\mathbb{P}^1)=p_2(L_2)$.
  
  Then there exists a normalized one-parameter subgroup $\lambda: G_m\rightarrow G$, where $\lim\limits_{t \to 0} \lambda(t).f=0$.
 $$\lambda(t)=\text{diag}(t^{r_0},t^{r_1})\times \text{diag} (t^{s_0},t^{s_1},t^{s_2}),$$
 such that $r_0+r_1=0$, $r_0\leq r_1$;   $s_0+s_1+s_2=0$, $s_0\leq s_1\leq s_2$; and $r_0+s_0<0$. 
 $$\lambda(t). f=\sum a_{\alpha\beta}(t^{r_0}x_0)^{\alpha_0}(t^{r_1}x_1)^{\alpha_1}(t^{s_0}y_0)^{\beta_0}(t^{s_1}y_1)^{\beta_1}(t^{s_2}y_2)^{\beta_2}=\sum a_{\alpha\beta}t^{r\alpha+s\beta}x^{\alpha}y^{\beta},$$
 $\lim\limits_{t \to 0} \lambda(t).f=0$, i.e., $r\alpha+s\beta>0$ for all $a_{\alpha\beta}\neq 0$. 

 It follows from the set of equations (\ref{equ3.2}) that if $S$ is unstable and has singularity at $P=[1,0]\times[1,0,0]$ which is not $A_1$-type, then the tangent cone of $S$ at $P$  is  $a_{22}y_2^2+b_{02}x_1y_2+c_{00}x_1^2$,  $a_{22}y_2^2+a_{12}y_1y_2+b_{02}x_1y_2$, or $a_{11}y_1^2+a_{22}y_2^2+a_{12}y_1y_2$.
 We consider these cases separately. 
 
 \bf Case-I: \rm The tangent cone of $S$ at $P$ is $f(1,x_1,1,y_1,y_2)_2=a_{22}y_2^2+b_{02}x_1y_2+c_{00}x_1^2$ where $c_{00}\neq0$. Note that $a_{22}\neq 0$, otherwise it contradicts the irreducibility of $f$. $c_{00}\neq0$ and $a_{22}\neq 0$ imply the inequality $s_0+s_2>0$. If $b_{11}\neq0$ then $s_1>0$ which contradicts $s_0+s_2>0$. Hence $b_{11}=0$ and 

$f(x_0,x_1,y_0,y_1,y_2)=a_{22}x_0^2y_2^2+x_0x_1(b_{22}y_2^2+b_{02}y_0y_2+b_{12}y_1y_2)+x_1^2(c_{00}y_0^2+c_{11}y_1^2+c_{22}y_2^2+c_{01}y_0y_1+c_{02}y_0y_2+c_{12}y_1y_2).$

The fibre over $[1,0]$ is the non-reduced conic $[1,0]\times Z(y_2^2)$. Note that  the map $[1,0]\times Z(y_2)=p_2^{-1}(Z(y_2))\rightarrow Z(y_2)$ is an isomorphism in some open subset of $Z(y_2)$. Therefore the generic points of the fibre over $P_1$ are in the ramification locus of $p_2$. 

\bf Case-II: \rm The tangent cone of $S$ at $P$ is $f(1,x_1,1,y_1,y_2)_2=a_{22}y_2^2+a_{12}y_1y_2+b_{02}x_1y_2$, where both $a_{12}$  and   $b_{02}$ are non-zero. These imply  inequalities
\begin{align}\label{equ2.4}
 2r_0+s_1+s_2>0,
\end{align}
\begin{align}\label{equ2.5}
 s_0+s_2>0.
\end{align}
But $b_{11}\neq0$ contradicts the inequality (\ref{equ2.5}) and $c_{01}\neq0 $ implies the inequality $2r_1+s_0+s_1>0$. This and the inequality (\ref{equ2.4}) imply $s_1>0$ which contradicts the inequality (\ref{equ2.5}). Hence $b_{11}=c_{01}=0$ and 

$f(x_0,x_1,y_0,y_1,y_2)=x_0^2(a_{22}y_2^2+a_{12}y_1y_2)+x_0x_1(b_{22}y_2^2+b_{02}y_0y_2+b_{12}y_1y_2)+x_1^2(c_{11}y_1^2+c_{22}y_2^2+c_{02}y_0y_2+c_{12}y_1y_2).$

Here  the generic points of the fibre component $[1,0]\times Z(y_2)$ over $[1,0]$ are in the ramification locus of $p_2$ and $p_2$ is not finite at $[1,0]\times[1,0,0]$. Note that there is a section $\sigma:\mathbb{P}^1\rightarrow Z(y_1,y_2) \simeq\mathbb{P}^1\times [1,0,0]\subset S $ of $p_1$. Moreover the corresponding map $\phi_{\sigma}:\mathbb{P}^1\rightarrow \mathbb{P}(T_{\mathbb{P}^2,P_2})$ is constant and for any $[\alpha_0,\alpha_1] \in {\mathbb P}^1$, $\phi_{\sigma}([\alpha_0,\alpha_1])=Z(y_2)$.

\bf Case-III: \rm The tangent cone of $S$ at $P$ is  $f(1,x_1,1,y_1,y_2)_2=a_{11}y_1^2+a_{22}y_2^2+a_{12}y_1y_2$. Then the tangent cone of $S$ at $P$ is a pull back of a tangent cone at $P_2$ by the map $p_2$.
\end{proof}
\end{lem}
Now we give a complete description of irreducible semi-stable surfaces in the following theorem.

\subsection{Proof of Theorem \ref{thm1.1}}
  Let $S\in \vert(2,2)\vert$ of $\mathbb{P}^1\times\mathbb{P}^2$ and $S=Z(f)$.
   We prove this result by contradiction. If
  $S$ is an irreducible unstable surface, then we proved in Lemma \ref{lem2.19}, that  some  $P=P_1\times P_2\in S_{\text{sing}}$  satisfy  one of the conditions (i), (ii), (iii) in the proof of Lemma \ref{lem2.19}.
  
 Conversely, we prove that if some singular points $P=P_1\times P_2\in S=Z(f)$ satisfy one of these conditions, then $S$ is unstable.
After a possible coordinate change, assume that $S$ passes through $P=[1,0]\times [1,0,0]$ and $P\in S_{\text{sing}}$. Then \\
$f(x_0,x_1,y_0,y_1,y_2)=x_0^2(a_{11}y_1^2+a_{22}y_2^2+a_{12}y_1y_2)+x_0x_1(b_{11}y_1^2+b_{22}y_2^2+b_{01}y_0y_1+b_{02}y_0y_2+b_{12}y_1y_2)+x_1^2(c_{00}y_0^2+c_{11}y_1^2+c_{22}y_2^2+c_{01}y_0y_1+c_{02}y_0y_2+c_{12}y_1y_2).$

If $f$ satisfies the condition (i) i.e., the tangent cone at $P$ is a pull back of a tangent cone at $P_2$, then 
$$f(1,0,1,y_1,y_2)_2=f(1,x_1,1,y_1,y_2)_2=a_{11}y_1^2+a_{22}y_2^2+a_{12}y_1y_2.$$
Therefore $f(x_0,x_1,y_0,y_1,y_2)=x_0^2(a_{11}y_1^2+a_{22}y_2^2+a_{12}y_1y_2)+x_0x_1(b_{11}y_1^2+b_{22}y_2^2+b_{12}y_1y_2)+x_1^2(c_{11}y_1^2+c_{22}y_2^2+c_{01}y_0y_1+c_{02}y_0y_2+c_{12}y_1y_2)$.
Let $\lambda(t)=\text{diag}(t^{-4},t^4)\times \text{diag}(t^{-10},t^5,t^5)$. Then clearly,
$$\lim\limits_{t\rightarrow 0}\lambda(t).f=0.$$
Hence $f$ is unstable.

Now assume that the condition (ii) holds, i.e.,  the fibre over $P_1\in\mathbb{P}^1$ is non-reduced and generic points of $p_1^{-1}(P_1)$ are in the ramification locus of $p_2$. Without loss of generality we may assume that the non-reduced fibre over $P_1=[1,0]$ is the conic $Z(y_2^2)$ and generically it is contained in the ramification locus of the map $p_2$. Therefore $b_{11}=b_{01}=0$. Then  \\
$f(x_0,x_1,y_0,y_1,y_2)=a_{22}x_0^2y_2^2+x_0x_1(b_{22}y_2^2+b_{02}y_0y_2+b_{12}y_1y_2)+x_1^2(c_{00}y_0^2+c_{11}y_1^2+c_{22}y_2^2+c_{01}y_0y_1+c_{02}y_0y_2+c_{12}y_1y_2).$ Let $\lambda(t)= \text{diag}(t^{-3},t^3)\times \text{diag}(t^{-2},t^{-2},t^{4})$. Then
$$\lim\limits_{t\rightarrow 0}\lambda(t).f=0.$$
Hence $S$ is an unstable surface.

If the condition (iii) holds, i.e., the fibre over $P_1$ is reduced and non-irreducible and  one of the components of the fibre over $P_1$ say $L_2$ is generically contained in the ramification locus of $p_2$, then\\
$f(x_0,x_1,y_0,y_1,y_2)=x_0^2(a_{22}y_2^2+a_{12}y_1y_2)+x_0x_1(b_{22}y_2^2+b_{02}y_0y_2+b_{12}y_1y_2)+x_1^2(c_{00}y_0^2+c_{11}y_1^2+c_{22}y_2^2+c_{01}y_0y_1+c_{02}y_0y_2+c_{12}y_1y_2).$\\
 Moreover $p_2$ is not finite at $P_1\times P_2=[1,0]\times [1,0,0]\in L_2$. Then   there is a section $\sigma:\mathbb{P}^1\rightarrow Z(y_1,y_2)\subset S$ of $p_1$  such that the corresponding map $\phi_{\sigma}:\mathbb{P}^1\rightarrow \mathbb{P}(T_{\mathbb{P}^2,P_2})$ is constant and  $\phi_{\sigma}(\mathbb{P}^1)=p_2(Z(x_1,y_2))=z(y_2)$. Hence $c_{00}=c_{01}=0$.

Therefore, $f(x_0,x_1,y_0,y_1,y_2)=x_0^2(a_{22}y_2^2+a_{12}y_1y_2)+x_0x_1(b_{22}y_2^2+b_{02}y_0y_2+b_{12}y_1y_2)+x_1^2(c_{11}y_1^2+c_{22}y_2^2+c_{02}y_0y_2+c_{12}y_1y_2)$.
Let $\lambda(t)=\text{diag}(t^{-2},t^2)\times \text{diag}(t^{-5},t^{-1},t^6)$. Then
$$\lim\limits_{t\rightarrow 0}\lambda(t).f=0.$$
Therefore $f$ is unstable.
 Hence  the theorem is proved. $\hfill\square$

\section{Irreducible stable surfaces}
In this section we describe irreducible stable surfaces. We know from \cite{MFK}, Proposition 4.2 that smooth effective divisors of $\mathbb{P}^n$ are stable. But here in our case not all smooth surfaces are stable. In the following lemma we describe when surfaces with at most $A_1$-type singularities are stable.
\begin{lem}\label{stablem2.18}
 Let $S=Z(f)$ be a $(2,2)$-type  surface  of $\mathbb{P}^1\times\mathbb{P}^2$ having at most $A_1$-type singularities. Then $S$ is stable if and only if one of the following conditions holds:
 
 (i) The second projection $p_2:S\rightarrow \mathbb{P}^2$ is a finite map.
 
 (ii) If $p_2$ is not a finite map, then it contracts some sections $\sigma:\mathbb{P}^1\rightarrow S$ of the first projection $p_1$ to $P_2\in\mathbb{P}^2$. Then $\phi_{\sigma}:\mathbb{P}^1\rightarrow \mathbb{P}(T_{\mathbb{P}^2, P_2})$ (see Notations and conventions \ref{item7} for the description of the map $\phi_{\sigma}$) is a non-constant map for all such sections.
 \begin{proof}
  We prove that if $S=Z(f)$ is not stable, then $p_2$ is not a finite map at $P=P_1\times P_2 \in S$. Hence there exists  a section $\sigma:\mathbb{P}^1\rightarrow \mathbb{P}^1\times P_2\subset S$, $p_1\circ \sigma=\text{id}$. Then it follows that   $p_2(p_1^{-1}(P_1))=C_{P_1}$ is a  conic which passes through $P_2$. Moreover, we prove that the morphism $\phi_{\sigma}:\mathbb{P}^1\rightarrow \mathbb{P}(T_{\mathbb{P}^2, P_2})$ is constant.
  
 Let $S$ be an irreducible non-stable surface having at most $A_1$-type singularities. Then there exists a normalized one-parameter subgroup $\lambda$ such that $\lim\limits_{t\rightarrow 0} \lambda(t).f$ exists. 
 $$\lambda(t).f=\sum a_{\alpha\beta}(t^{r_0}x_0)^{\alpha_0}(t^{r_1}x_1)^{\alpha_1}(t^{s_0}y_0)^{\beta_0}(t^{s_1}y_1)^{\beta_1}(t^{s_2}y_2)^{\beta_2}=\sum a_{\alpha\beta}t^{r\alpha+s\beta}x^{\alpha}y^{\beta},$$
 then $r\alpha+s\beta\geq0$ for all $a_{\alpha\beta}\neq 0$. Note that if $s_0=s_1=s_2=0$, then $r_0<r_1$. Then $r\alpha+s\beta=r_0\alpha_0+r_1\alpha_1\geq0$ if and only if $\alpha_1\neq0$ for all $a_{\alpha\beta}\neq 0$. Therefore, $f$ is divisible by $x_1$. But this contradicts the irreducibility of $f$. Hence $s_0=s_1=s_2=0$ is not possible as $S$ is not stable.

An explicit form of the (2,2)-type polynomial $f$ is 
 
 $f(x_0,x_1,y_0,y_1,y_2)=x_0^2(a_{11}y_1^2+a_{22}y_2^2+a_{12}y_1y_2+a_{01}y_0y_1+a_{02}y_0y_2+a_{00}y_0^2)+x_0x_1(b_{11}y_1^2+b_{22}y_2^2+b_{01}y_0y_1+b_{02}y_0y_2+b_{12}y_1y_2+b_{00}y_0^2)+x_1^2(c_{00}y_0^2+c_{11}y_1^2+c_{22}y_2^2+c_{01}y_0y_1+c_{02}y_0y_2+c_{12}y_1y_2).$
 
 As $S$ is non-stable, the coefficient of $x_0^2y_0^2$, $a_{00}=0$, otherwise $r_0+s_0\geq0$ which is a contradiction.
 The coefficient of $x_0^2y_0y_1$, $a_{01}=0$, otherwise $2r_0+s_0+s_1\geq 0$ which implies $s_0+s_1\geq0$. But we know $s_2\geq 0$. Hence there is a contradiction .
 Similarly, the coefficient of $x_0x_1y_0^2$, $b_{00}=0$, otherwise $r_0+r_1+2s_0\geq0$ which implies $s_0\geq0$. But this gives
   a contradiction.
Also, the coefficient of $x_0x_1y_0y_1$, $b_{01}=0$, otherwise $r_0+r_1+s_0+s_1\geq 0$ which implies $s_2\leq 0$. But this is a contradiction. 

Now we rewrite the polynomial $f$,
 
 $f(x_0,x_1,y_0,y_1,y_2)=x_0^2(a_{11}y_1^2+a_{22}y_2^2+a_{12}y_1y_2+a_{02}y_0y_2)+x_0x_1(b_{11}y_1^2+b_{22}y_2^2+b_{02}y_0y_2+b_{12}y_1y_2)+x_1^2(c_{00}y_0^2+c_{11}y_1^2+c_{22}y_2^2+c_{01}y_0y_1+c_{02}y_0y_2+c_{12}y_1y_2)$.
 
  If $S$ is smooth at $P=[1,0]\times [1,0,0]$, then $a_{02}\neq 0$. Then the non-stability of $f$ implies the inequality
  \begin{align}\label{equ2.3}
   2r_0+s_0+s_2\geq0.
  \end{align}
  If the coefficient of $x_1^2y_0^2$, $c_{00}\neq 0$, then $2r_1+2s_0\geq0$. This inequality and the inequality  (\ref{equ2.3}) imply $3s_0+s_2\geq0$. But note that $s_0+s_1+s_2=0$. Hence  $2s_0-s_1\geq0$ implies $2s_0\geq s_1\geq s_0$ and $s_0=s_1=s_2=0$, which is a contradiction.
  
  If the coefficient of $x_1^2y_0y_1$, $c_{01}\neq 0$, then $2r_1+s_0+s_1\geq0$. This and inequality (\ref{equ2.3}) imply $s_0\geq0$ which is a contradiction.
 
If $S$ is smooth but not stable, then $c_{00}=c_{01}=0$. Hence,
$f(x_0,x_1,y_0,y_1,y_2)=x_0^2(a_{11}y_1^2+a_{22}y_2^2+a_{12}y_1y_2+a_{02}y_0y_2)+x_0x_1(b_{11}y_1^2+b_{22}y_2^2+b_{02}y_0y_2+b_{12}y_1y_2
)+x_1^2(c_{11}y_1^2+c_{22}y_2^2+c_{02}y_0y_2+c_{12}y_1y_2)$.

Note that $p_2$ is not finite at the point $[1,0]\times [1,0,0]$. Then this corresponds to the section $\sigma:\mathbb{P}^1\rightarrow Z(y_1,y_2)\simeq\mathbb{P}^1\times [1,0,0]\subset S$  of $p_1$. Let 

$\varphi=(\alpha_0^2a_{11}+\alpha_0\alpha_1b_{11}+\alpha_1^2c_{11})y_1^2+(\alpha_0^2a_{22}+\alpha_0\alpha_1b_{22}+\alpha_1^2c_{22})y_2^2+(\alpha_0^2a_{12}+\alpha_0\alpha_1b_{12}+\alpha_1^2c_{12})y_1y_2+(\alpha_0^2a_{02}+\alpha_0\alpha_1b_{02}+\alpha_1^2c_{02})y_0y_2$.

$p_2(p_1^{-1}([\alpha_0,\alpha_1]))=p_2([\alpha_0,\alpha_1]\times Z(\varphi))=Z(\varphi)$ is the  conic passing through $[1,0,0]$. Note that the coefficient of $y_0y_2$, $\alpha_0^2a_{02}+\alpha_0\alpha_1b_{12}+\alpha_2^2c_{02}\neq0$, otherwise $S$ becomes singular at $[\alpha_0,\alpha_1]\times [1,0,0]$. Then $Z(y_2)$ is the tangent line of all conic $p_2([\alpha_0,\alpha_1]\times Z(\varphi))=Z(\varphi)$ at $[1,0,0]$. Hence the map $\phi_{\sigma}:\mathbb{P}^1\rightarrow \mathbb{P}(T_{\mathbb{P}^2, P_2})$ is constant. 

If $S$  has a singularity at $P=[1,0]\times[1,0,0]$, then $a_{02}=0$, i.e.,  $f(1,x_1,1,y_1,y_2)_1=0$. Then the tangent cone of $S$ at $P$
 is $$f(1,x_1,1,y_1,y_2)_2=a_{11}y_1^2+a_{22}y_2^2+a_{12}y_1y_2+b_{02}x_1y_2+c_{00}x_1^2.$$ 
 Note that, from the hypothesis $S$ has an  $A_1$-type singularity at $P$.
So the determinant of the Hessian of $f$ at $P$, $H=8a_{11}a_{22}c_{00}-2a_{11}b_{02}^2-2a_{12}^2c_{00}\neq 0$.

If $a_{11}\neq 0$ and $c_{00}\neq 0$ then $2r_0+2s_1\geq0$ and $2r_1+2s_0\geq0$ respectively. They imply $s_2\leq0$, which is a contradiction. So $a_{11}$ and $c_{00}$ are not simultaneously non-zero.

If $a_{11}=0$, then $f(1,x_1,1,y_1,y_2)_2=a_{22}y_2^2+a_{12}y_1y_2+b_{02}x_1y_2+c_{00}x_1^2$ and $H=-2a_{12}^2c_{00}.$

If $c_{00}=0$, then
$f(1,x_1,1,y_1,y_2)_2=a_{11}y_1^2+a_{22}y_2^2+a_{12}y_1y_2+b_{02}x_1y_2$ and $H=-2a_{11}b_{02}^2$.

If $a_{11}=0$, then $a_{12}\neq0$ and $c_{00}\neq0$ as $H\neq0$. Hence we have following inequalities
$2r_0+s_1+s_2\geq0$, and  $2r_1+2s_0\geq0$ respectively. They imply $r_0\geq0$, i.e., $r_0=r_1=0$ and $s_0\geq0$ which is a contradiction.

Similarly, if $c_{00}=0$, then $a_{11}\neq0$ and $b_{02}\neq0$  as $H\neq0$. Hence we have following inequalities
$r_0+s_1\geq0$ and  $r_0+r_1+s_0+s_2\geq0$. They imply $r_0\geq0$, hence $r_0=r_1=0$ and $s_1=s_0+s_2=0$. If the coefficient of $x_1^2y_0y_1$, $c_{01}\neq0$ then $2r_1+s_0+s_1\geq0$ which is a contradiction. Then $c_{01}=0$ and finally $f(x_0,x_1,y_0,y_1,y_2)=x_0^2(a_{11}y_1^2+a_{22}y_2^2+a_{12}y_1y_2)+x_0x_1(b_{11}y_1^2+b_{22}y_2^2+b_{02}y_0y_2+b_{12}y_1y_2)+x_1^2(c_{11}y_1^2+c_{22}y_2^2+c_{02}y_0y_2+c_{12}y_1y_2)$.

$L'=\mathbb{P}^1\times [1,0,0]\subset S$  corresponds to a section $\sigma$ of $p_1$ and $p_2(L')= [1,0,0]$. Let

$\varphi=(\alpha_0^2a_{11}+\alpha_0\alpha_1b_{11}+\alpha_1^2c_{11})y_1^2+(\alpha_0^2a_{22}+\alpha_0\alpha_1b_{22}+\alpha_1^2c_{22})y_2^2+(\alpha_0^2a_{12}+\alpha_0\alpha_1b_{12}+\alpha_1^2c_{12})y_1y_2+(\alpha_0\alpha_1b_{02}+\alpha_1^2c_{02})y_0y_2$.

Therefore, $p_2(p_1^{-1}([\alpha_0,\alpha_1]))= Z(\varphi)$. Except finitely many points in $\mathbb{P}^1$, the conic $p_2(p_1^{-1}([\alpha_0,\alpha_1]))$ is smooth at  $P_2=[1,0,0]$ and $Z(y_2)$ is the tangent of  the conic  $Z(\varphi)$ at $P_2=[1,0,0]$. Hence the map $\phi_{\sigma}:\mathbb{P}^1\rightarrow \mathbb{P}(T_{\mathbb{P}^2, P_2})$ is constant.

Now we prove the converse direction. Let us assume that  $p_2:S\rightarrow \mathbb{P}^2$ is not a finite map.
After some coordinate change we may assume that $S$ passes through $P=[1,0]\times[1,0,0]$ and $p_2$ is not finite at $P$. Then there is  a section $\sigma:\mathbb{P}^1\rightarrow Z(y_1,y_2)\subset S$ of $p_1$ which passes through $P$. Hence

$f(x_0,x_1,y_0,y_1,y_2)=x_0^2(a_{11}y_1^2+a_{22}y_2^2+a_{12}y_1y_2+a_{01}y_0y_1+a_{02}y_0y_2)+x_0x_1(b_{11}y_1^2+b_{22}y_2^2+b_{01}y_0y_1+b_{02}y_0y_2+b_{12}y_1y_2)+x_1^2(c_{11}y_1^2+c_{22}y_2^2+c_{01}y_0y_1+c_{02}y_0y_2+c_{12}y_1y_2).$

We prove that if the map $\phi_{\sigma}:\mathbb{P}^1\rightarrow \mathbb{P}(T_{\mathbb{P}^2, P_2})$ is constant then $S$ is not a stable surface. 

Now consider $p_2(p_1^{-1}[\alpha_0,\alpha_1])=Z(\varphi)=Z((\alpha_0^2a_{11}+\alpha_0\alpha_1b_{11}+\alpha_1^2c_{11})y_1^2+(\alpha_0^2a_{22}+\alpha_0\alpha_1b_{22}+\alpha_1^2c_{22})y_2^2+(\alpha_0^2a_{12}+\alpha_0\alpha_1b_{12}+\alpha_1^2c_{12})y_1y_2+(\alpha_0^2a_{02}+\alpha_0\alpha_1b_{02}+\alpha_1^2c_{02})y_0y_2+(\alpha_0^2a_{01}+\alpha_0\alpha_1b_{01}+\alpha_1^2c_{01})y_0y_1)$ which is a  conic passing through $[1,0,0]$ and generic conics are smooth. Therefore $Z((\alpha_0^2a_{02}+\alpha_0\alpha_1b_{02}+\alpha_1^2c_{02})y_2+(\alpha_0^2a_{01}+\alpha_0\alpha_1b_{01}+\alpha_1^2c_{01})y_1)$ is the tangent line of smooth conics at the   point $[1,0,0]$. As $\phi_{\sigma}([\alpha_0,\alpha_1])=Z((\alpha_0^2a_{02}+\alpha_0\alpha_1b_{02}+\alpha_1^2c_{02})y_2+(\alpha_0^2a_{01}+\alpha_0\alpha_1b_{01}+\alpha_1^2c_{01})y_1)$ is a constant map,  without loss of generality we may assume that that $\phi_{\sigma}([\alpha_0,\alpha_1])=Z(y_2)$.

Therefore $a_{01}=c_{01}=b_{01}=0$, and

$f(x_0,x_1,y_0,y_1,y_2)=x_0^2(a_{11}y_1^2+a_{22}y_2^2+a_{12}y_1y_2+a_{02}y_0y_2)+x_0x_1(b_{11}y_1^2+b_{22}y_2^2+b_{02}y_0y_2+b_{12}y_1y_2)+x_1^2(c_{11}y_1^2+c_{22}y_2^2+c_{02}y_0y_2+c_{12}y_1y_2).$

Let $\lambda(t)=\text{diag}(t^{0},t^{0})\times \text{diag} (t^{-1},t^{0},t^{1}).$ 
 Then $\lim\limits_{t\rightarrow 0} \lambda(t).f$ exists. Hence $S$ is a non-stable surface.
\end{proof} 
\end{lem}
Now we describe necessary and sufficient conditions of a $(2,2)$-type semi-stable surface  to be a stable surface.

\subsection{Proof of Theorem \ref{thm1.2}}
Let  $S=Z(f)$ be a semi-stable irreducible surface.
We prove that if $S$ is not stable, then one of the following conditions holds: 
\\
(i) The map $p_2$ is not finite at  $P=P_1\times P_2\in S$, where $P$ is either smooth or $A_1$-type singular point of $S$. Then there is a section $\sigma:\mathbb{P}^1\rightarrow S$ of $p_1$ contracted to $P_2$. Moreover the corresponding map $\phi_{\sigma}:\mathbb{P}^1\rightarrow \mathbb{P}(T_{\mathbb{P}^2,P_2})$ is constant.\\
  (ii) The map $p_2$ is not finite at $P\in S_{\text{sing}}$ where $P$ is not $A_1$-type singular point of $S$. \\
  (iii) The fibre over $P_1$ is non-reduced where $P=P_1\times P_2\in S_{\text{sing}}$ and $P$ is not $A_1$-type singular point of $S$.\\
 Let $\lambda$ be a normalized one-parameter subgroup such that $\lim\limits_{t\rightarrow 0} \lambda(t).f$ exists.
 
 We proved in Lemma \ref{stablem2.18} that if $f$ is non-stable then the coefficients of $x_0^2y_0^2$, $x_0^2y_0y_1$, $x_0x_1y_0^2$ and $x_0x_1y_0y_1$ are zero. Then 
 
 $f(x_0,x_1,y_0,y_1,y_2)=x_0^2(a_{11}y_1^2+a_{22}y_2^2+a_{12}y_1y_2+a_{02}y_0y_2)+x_0x_1(b_{11}y_1^2+b_{22}y_2^2+b_{02}y_0y_2+b_{12}y_1y_2)+x_1^2(c_{00}y_0^2+c_{11}y_1^2+c_{22}y_2^2+c_{01}y_0y_1+c_{02}y_0y_2+c_{12}y_1y_2).$
 
 If $S$ is either smooth or has singularity of type $A_1$ at $P=[1,0]\times [1,0,0]$, then we proved  in Lemma \ref{stablem2.18} that $p_2$ is not finite at  $P$. Then there is a section $\sigma:\mathbb{P}^1\rightarrow \mathbb{P}^1\times [1,0,0]\subset S$ of $p_1$ contracted to $P_2=[1,0,0]$. Moreover the corresponding map $\phi_{\sigma}:\mathbb{P}^1\rightarrow \mathbb{P}(T_{\mathbb{P}^2,P_2})$ is constant. Now consider that $P$ is a singular point of $S$ which is not $A_1$-type. Hence,
  $f(1,x_1,1,y_1,y_2)_1=0$ and 
 the tangent cone of $S$ at $P$ is $$f(1,x_1,1,y_1,y_2)_2=a_{11}y_1^2+a_{22}y_2^2+a_{12}y_1y_2+b_{02}x_1y_2+c_{00}x_1^2$$ 
  and the determinant of the Hessian of $f$ at $P$ is $H=8a_{11}a_{22}c_{00}-2a_{11}b_{02}^2-2a_{12}^2c_{00}=0$. 

If $a_{11}\neq 0$ and $c_{00}\neq 0$ then $2r_0+2s_1\geq0$ and $2r_1+2s_0\geq0$ respectively. They imply $s_2\leq0$ which is a contradiction. Therefore either $a_{11}=0$ or $c_{00}=0$.

If $a_{11}=0$, then $f(1,x_1,1,y_1,y_2)_2=a_{22}y_2^2+a_{12}y_1y_2+b_{02}x_1y_2+c_{00}x_1^2$ and $H=-2a_{12}^2c_{00}=0$. Then either $a_{12}=0$ or $c_{00}=0$.

If $c_{00}=0$, then
$f(1,x_1,1,y_1,y_2)_2=a_{11}y_1^2+a_{22}y_2^2+a_{12}y_1y_2+b_{02}x_1y_2$ and $H=-2a_{11}b_{02}^2=0$. This implies $a_{11}=0$. Note that $b_{02}\neq0$, otherwise $f(1,0,1,y_1,y_2)_2=f(1,x_1,1,y_1,y_2)_2=a_{11}y_1^2+a_{22}y_2^2+a_{12}y_1y_2$. Therefore $S$ becomes non-stable (see Theorem \ref{thm1.1}(i)). Therefore, either $f(x_0,x_1,y_0,y_1,y_2)=x_0^2(a_{22}y_2^2+a_{12}y_1y_2)+x_0x_1(b_{11}y_1^2+b_{22}y_2^2+b_{02}y_0y_2+b_{12}y_1y_2)+x_1^2(c_{11}y_1^2+c_{22}y_2^2+c_{01}y_0y_1+c_{02}y_0y_2+c_{12}y_1y_2)$ or $f(x_0,x_1,y_0,y_1,y_2)=a_{22}x_0^2y_2^2+x_0x_1(b_{11}y_1^2+b_{22}y_2^2+b_{02}y_0y_2+b_{12}y_1y_2)+x_1^2(c_{00}y_0^2+c_{11}y_1^2+c_{22}y_2^2+c_{01}y_0y_1+c_{02}y_0y_2+c_{12}y_1y_2)$.

In the former case, the section $\sigma:\mathbb{P}^1\rightarrow Z(y_1,y_2)$ of $p_1$ passes through $P$. Hence the condition (ii) holds.

In the latter case, the fibre over $[1,0]$ is non-reduced. Hence the condition (iii) holds.

Now we prove the converse direction of the theorem which is if $S$ satisfies (i), (ii) or (iii), then $S$ is  non-stable. 
  In Lemma \ref{stablem2.18}, we proved that if the condition (i) holds, then $S$ is non-stable.
Now let $P=P_1\times P_2$ be a singular point of $S$ which is not $A_1$-type.
After a suitable coordinate change we may assume  that  $P=[1,0]\times[1,0,0]$. Then 

$f(x_0,x_1,y_0,y_1,y_2)=x_0^2(a_{11}y_1^2+a_{22}y_2^2+a_{12}y_1y_2)+x_0x_1(b_{11}y_1^2+b_{22}y_2^2+b_{02}y_0y_2+b_{12}y_1y_2)+x_1^2(c_{00}y_0^2+c_{11}y_1^2+c_{22}y_2^2+c_{01}y_0y_1+c_{02}y_0y_2+c_{12}y_1y_2).$
 The tangent cone of $S$ at $P$ is
$$f(1,x_1,1,y_1,y_2)_2=a_{11}y_1^2+a_{22}y_2^2+a_{12}y_1y_2+b_{02}x_1y_2+c_{00}x_1^2,$$ 
and the determinant of the Hessian matrix at the point $P$ is $H=8a_{11}a_{22}c_{00}-2a_{11}b_{02}^2-2a_{12}^2c_{00}=0$.

If $p_2$ is not finite at $P=[1,0]\times [1,0,0]$, then  the section $\sigma:\mathbb{P}^1\rightarrow Z(y_1,y_2)$ passes through $P$. Therefore $c_{00}=0$ and $H=-2a_{11}b_{02}^2=0$. As $S$ is semi-stable, $b_{02}\neq0$. Then $a_{11}=0$. Now  $f(1,x_1,1,y_1,y_2)_2=a_{22}y_2^2+a_{12}y_1y_2+b_{02}x_1y_2$, and  $f(x_0,x_1,y_0,y_1,y_2)=x_0^2(a_{22}y_2^2+a_{12}y_1y_2)+x_0x_1(b_{11}y_1^2+b_{22}y_2^2+b_{02}y_0y_2+b_{12}y_1y_2)+x_1^2(c_{11}y_1^2+c_{22}y_2^2+c_{01}y_0y_1+c_{02}y_0y_2+c_{12}y_1y_2).$
Let $\lambda(t)=\text{diag}(t^{-1},t)\times \text{diag}(t^{-2},t^0,t^2)$. Then clearly,
$\lim\limits_{t \to 0}\lambda(t).f$ exists. Hence $f$ is not stable.

If the fibre over $[1,0]$ is non-reduced, then without loss of generality we may assume that \\
$f(x_0,x_1,y_0,y_1,y_2)=a_{22}x_0^2y_2^2+x_0x_1(b_{11}y_1^2+b_{22}y_2^2+b_{01}y_0y_1+b_{02}y_0y_2+b_{12}y_1y_2)+x_1^2(c_{00}y_0^2+c_{11}y_1^2+c_{22}y_2^2+c_{01}y_0y_1+c_{02}y_0y_2+c_{12}y_1y_2).$ Then $f(1,x_1,1,y_1,y_2)_2=a_{22}y_2^2+b_{01}y_0y_1+b_{02}x_1y_2+c_{00}x_1^2$. The Hessian at $P$ is $H=2a_{22}b_{01}^2=0$. As $a_{22}\neq0$, $b_{01}=0$. 

Now $f(1,x_1,1,y_1,y_2)_2=a_{22}y_2^2+b_{02}x_1y_2+c_{00}x_1^2$, and  $f(x_0,x_1,y_0,y_1,y_2)=a_{22}x_0^2y_2^2+x_0x_1(b_{11}y_1^2+b_{22}y_2^2+b_{02}y_0y_2+b_{12}y_1y_2)+x_1^2(c_{00}y_0^2+c_{11}y_1^2+c_{22}y_2^2+c_{01}y_0y_1+c_{02}y_0y_2+c_{12}y_1y_2).$
Let $\lambda(t)=\text{diag}(t^{-1},t)\times \text{diag}(t^{-1},t^0,t)$. Then clearly 
$\lim\limits_{t \to 0}\lambda(t).f$ exists. Hence $f$ is not stable. \\Therefore we proved the theorem.

In the next proposition, we see that there are some irreducible stable (2,2)-type surfaces which have higher order singularities.
 \begin{prop}\label{prop3.9}
 Let $S=Z(f)$ be an irreducible, singular, stable $(2,2)$-type surface and $S$ has quasi-homogeneous germs (Definition \ref{defi2.1}) at singular points. Then $S$ has either $A_1$, $A_2$ or $A_3$-type isolated singularities or non-isolated singularities. 
\begin{proof}  
 If $f(x_0,x_1,y_0,y_1,y_2)=x_0^2(a_{11}y_1^2+a_{12}y_1y_2)+x_0x_1(b_{11}y_1^2+b_{22}y_2^2+b_{02}y_0y_2+b_{12}y_1y_2)+x_1^2(c_{00}y_0^2+c_{11}y_1^2+c_{22}y_2^2+c_{02}y_0y_2+c_{12}y_1y_2)$, then $S=Z(f)$ is a stable surface from Theorem \ref{thm1.2} where $a_{11}, b_{02}, a_{12}$ and $c_{00}$ are non-zero.
Note that the determinant of the Hessian matrix at the point $P=[1,0]\times[1,0,0]$ is $H=-2a_{11}b_{02}^2-2a_{12}^2c_{00}$. 
If $H\neq0$, then $S$ has $A_1$-type singularity at $P$.
If $H=0$, then $S$ has other than $A_1$-type singularity at $P=[1,0]\times[1,0,0]$. Consider the affine neighbourhood  $\{x_0=1,y_0=1\}$ of $P$ and 

$f(1,x_1,1,y_1,y_2)=f(1,x_1,1,y_1,y_2)_2+f(1,x_1,1,y_1,y_2)_3+f(1,x_1,1,y_1,y_2)_4$\\
$=(a_{11}y_1^2+a_{12}y_1y_2+b_{02}x_1y_2+c_{00}x_1^2)+(b_{11}x_1y_1^2+b_{22}x_1y_2^2+b_{12}x_1y_1y_2+c_{02}x_1^2y_2)+(c_{11}x_1^2y_1^2+c_{22}x_1^2y_2^2+c_{12}x_1^2y_1y_2).$

The partial derivatives of $f(1,x_1,1,y_1,y_2)$ with respect to $x_1$, $y_1$ and $y_2$ are
 
$f_{x_1}=b_{02}y_2+2c_{00}x_1+b_{11}y_1^2+b_{22}y_2^2+b_{12}y_1y_2+2c_{02}x_1y_2+2c_{11}x_1y_1^2+2c_{22}x_1y_2^2+2c_{12}x_1y_1y_2,$
  
$f_{y_1}=2a_{11}y_1+a_{12}y_2+2b_{11}x_1y_1+b_{12}x_1y_2+2c_{11}x_1^2y_1+c_{12}x_1^2y_2,$
  
$f_{y_2}=a_{12}y_1+b_{02}x_1+2b_{22}x_1y_2+b_{12}x_1y_1+c_{02}x_1^2+2c_{22}x_1^2y_2+c_{12}x_1^2y_1$.
   
As $a_{11}b_{02}^2+a_{12}^2c_{00}=0$, the linear terms of partial derivatives of $f$ are not linearly independent. Let
 
$g_1=2a_{12}c_{00}f_{y_2}+b_{02}^2f_{y_1}-b_{02}a_{12}f_{x_1}=x_1^2(2a_{12}c_{00}c_{02}+y_1(2a_{12}c_{00}c_{12}+2c_{00}b_{02}^2)+y_2(4a_{12}c_{00}c_{22}+b_{02}^2c_{12}))+y_1^2(-b_{11}b_{02}a_{12}-2a_{12}b_{02}c_{11}x_1)+y_2^2(-b_{22}b_{02}a_{12}-2a_{12}b_{02}c_{22}x_1)+x_1y_1(2a_{12}b_{12}c_{00}+2b_{11}b_{02}^2)+x_1y_2(4a_{12}b_{22}c_{00}+b_{12}b_{02}^2-2a_{12}b_{02}c_{02})+y_1y_2(-a_{12}b_{02}b_{12}-2a_{12}b_{02}c_{12}x_1)$, 
 
 $g_2=f_{y_1}=\beta_1y_1+\beta_2y_2$, $\beta_1=(2a_{11}+2b_{11}x_1+2c_{11}x_1^2)$ and $\beta_2=(a_{12}+b_{12}x_1+c_{12}x_1^2)$ are units in  $\mathbb{C}[[x_1,y_1,y_2]]$,
 
 $g_3=f_{y_2}=\gamma_1y_1+\gamma_2x_1$, $\gamma_1=(a_{12}+b_{12}x_1+c_{12}x_1^2)$ and $\gamma_2=(b_{02}+2b_{22}y_2+c_{02}x_1+2c_{22}x_1y_2)$ are units in  $\mathbb{C}[[x_1,y_1,y_2]]$.
 
 Note that  $S$ has a quasi-homogeneous germ at $P$ which is same as saying\\ $(Z(f(1,x_1,1,y_1,y_2),(0,0,0)))$ is quasi-homogeneous i.e.,\\
 $f(1,x_1,1,y_1,y_2)\in (f_{x_1},f_{y_1},f_{y_2})=(g_1,g_2,g_3)$.
 
 Let us consider a ring homomorphism $\phi(y_1)=y_1$, $\phi(g_2)=y_2$ and $\phi(g_3)=x_1$ of $\mathbb{C}[[x_1,y_1,y_2]]$. $\phi$ is an automorphism  from Theorem \ref{thm2.13}. $\alpha_1=-\frac{\beta_1}{\beta_2}$ and $\alpha_2=-\frac{\gamma_1}{\gamma_2}$ are  units in $\mathbb{C}[[x_1,y_1,y_2]]$. The following equation is the coefficient of the degree two part of the image of  $g_1$ in $\mathbb{C}[[x_1,y_1,y_2]]/(g_2,g_3)$.
 If 
  \begin{align}\label{equ3.30}
 \alpha_2^2(2a_{12}c_{00}c_{02})-b_{11}b_{02}a_{12}+\alpha_1^2(-b_{22}b_{02}a_{12})+\alpha_2(2a_{12}b_{12}c_{00}+2b_{11}b_{02}^2)+
 \end{align}
 \begin{align*}
 \alpha_1\alpha_2(4a_{12}b_{22}c_{00}+b_{12}b_{02}^2-2a_{12}b_{02}c_{02})-\alpha_1a_{12}b_{02}b_{12}\neq0,
\end{align*}
then the image of $g_1$ is $ay_1^2$ in $\mathbb{C}[[x_1,y_1,y_2]]/(g_2,g_3)$  where $a$ is a  unit in $\mathbb{C}[[x_1,y_1,y_2]]$. Hence
 $$\mathbb{C}[[x_1,y_1,y_2]]/(f_{x_1},f_{y_1},f_{y_2})\cong \mathbb{C}[[x_1,y_1,y_2]]/(g_1,g_2,g_3)\cong \mathbb{C}[[y_1]]/(y_1^2).$$
  Therefore $S$ has an $A_2$-type singularity at  point $P$ (see Proposition \ref{prop2.14}).

Note that following equation is the coefficient of the degree three part of image of  $g_1$  in $\mathbb{C}[[x_1,y_1,y_2]]/(g_2,g_3)$.
If equation (\ref{equ3.30}) is zero but 
\begin{align}\label{equ3.32}
 \alpha_2^2(2a_{12}c_{00}c_{12}+2c_{00}b_{02}^2)+\alpha_1\alpha_2^2(4a_{12}c_{00}c_{22}+b_{02}^2c_{12})-
\end{align}
\begin{align*}
\alpha_22a_{12}b_{02}c_{11}-\alpha_1^2\alpha_22a_{12}b_{02}c_{11}-\alpha_1\alpha_22a_{12}b_{02}c_{00}\neq0,
\end{align*} 
then image of $g_1$ is $by_1^3$ in $\mathbb{C}[[x_1,y_1,y_2]]/(g_2,g_3)$  where $b$ is a unit in $\mathbb{C}[[x_1,y_1,y_2]]$. Therefore
$$\mathbb{C}[[x_1,y_1,y_2]]/(f_{x_1},f_{y_1},f_{y_2})\cong \mathbb{C}[[x_1,y_1,y_2]]/(g_1,g_2,g_3)\cong \mathbb{C}[[y_1]]/(y_1^3).$$
Hence $S$ has an $A_3$ type singularity at $P$ (see Proposition \ref{prop2.14}).
Now assume equations (\ref{equ3.30}) and (\ref{equ3.32}) are zero. Then the image of $g_1$  is zero in $\mathbb{C}[[x_1,y_1,y_2]]/(g_2,g_3)$. Therefore
$$\mathbb{C}[[x_1,y_1,y_2]]/(f_{x_1},f_{y_1},f_{y_2})\cong \mathbb{C}[[x_1,y_1,y_2]]/(g_1,g_2,g_3)\cong \mathbb{C}[[y_1]].$$
Hence $S$ has a non-isolated singularity at $P$.
\end{proof}
 \end{prop}
\section{Irreducible strictly semi-stable surfaces}
Using Theorem \ref{thm1.2},  we  list all strictly semi-stable $(2,2)$-type surfaces and their degenerations in the following theorem.
\begin{thm}\label{thm3.11}
  Let $S=Z(f)$ be an irreducible semi-stable surface. Then $S$ is strictly semi-stable if and only if $S$ satisfies one of the following conditions:
  
  (i) There are some sections $\sigma:\mathbb{P}^1\rightarrow S$ of $p_1$ such that $\phi_{\sigma}:\mathbb{P}^1\rightarrow \mathbb{P}(T_{\mathbb{P}^2,P_2})$ (see Notations and conventions (\ref{item7}) for the description of $\phi_{\sigma}$) is a constant map, where $p_2(\sigma(\mathbb{P}^1))=P_2\in \mathbb{P}^2$.
  
  (ii) The map $p_2$ is not  finite at $P\in S_{\text{sing}}$, where $P$ is not an $A_1$-type singularity. 

(iii) $P=P_1\times P_2\in S_{\text{sing}}$ which is not an $A_1$-type and the fibre over $P_1$ is non-reduced.
\begin{proof}
Let $S$ be an irreducible, semi-stable surface and without loss of generality assume that $S$ passes through the point $P=[1,0]\times[1,0,0]$. Then $f(x_0,x_1,y_0,y_1,y_2)=x_0^2(a_{11}y_1^2+a_{22}y_2^2+a_{12}y_1y_2+a_{02}y_0y_2)+x_0x_1(b_{11}y_1^2+b_{22}y_2^2+b_{01}y_0y_1+b_{02}y_0y_2+b_{12}y_1y_2+b_{00}y_0^2)+x_1^2(c_{00}y_0^2+c_{11}y_1^2+c_{22}y_2^2+c_{01}y_0y_1+c_{02}y_0y_2+c_{12}y_1y_2).$

\bf (i) \rm Let $\sigma:\mathbb{P}^1\rightarrow Z(y_1,y_2)\subset S$ be a section of $p_1$ such that $\phi_{\sigma}:\mathbb{P}^1\rightarrow \mathbb{P}(T_{\mathbb{P}^2,P_2})$ is a constant map. Then from  Lemma \ref{stablem2.18},   $f=x_0^2(a_{11}y_1^2+a_{22}y_2^2+a_{12}y_1y_2+a_{02}y_0y_2)+x_0x_1(b_{11}y_1^2+b_{22}y_2^2+b_{02}y_0y_2+b_{12}y_1y_2)+x_1^2(c_{11}y_1^2+c_{22}y_2^2+c_{02}y_0y_2+c_{12}y_1y_2)$. 
 Note that $a_{11}\neq0$, $b_{02}\neq0$; and $c_{11}\neq0$ or $c_{02}\neq0$.  Otherwise it contradicts the semi-stability of $f$. Let  $\lambda(t)=\text{diag}(t^{0},t^0)\times \text{diag}(t^{-1},t^0,t)$. Then $\lim\limits_{t\rightarrow 0} \lambda(t).f=x_0^2(a_{11}y_1^2+a_{02}y_0y_2)+x_0x_1(b_{11}y_1^2+b_{02}y_0y_2)+x_1^2(c_{11}y_1^2+c_{02}y_0y_2)$.

\bf (ii) \rm Let $P$ be a singular point of $S$ which is not an $A_1$-type. Also assume that there is a section $\sigma:\mathbb{P}^1\rightarrow Z(y_1,y_2)\subset S$ of $p_1$ which passes through $P$. Hence from Theorem  \ref{thm1.2},  the polynomial $f(x_0,x_1,y_0,y_1,y_2)=x_0^2(a_{22}y_2^2+a_{12}y_1y_2)+x_0x_1(b_{11}y_1^2+b_{22}y_2^2+b_{02}y_0y_2+b_{12}y_1y_2)+x_1^2(c_{11}y_1^2+c_{22}y_2^2+c_{01}y_0y_1+c_{02}y_0y_2+c_{12}y_1y_2)$, where $a_{12}\neq0$, $b_{02}\neq0$ and $c_{01}\neq0$. Otherwise $S$ will become unstable. If we take the one-parameter subgroup $\lambda(t)=\text{diag}(t^{-1},t)\times \text{diag}(t^{-2},t^0,t^2)$, then $\lim\limits_{t\rightarrow 0} \lambda(t).f=a_{12}x_0^2y_1y_2+x_0x_1(b_{11}y_1^2+b_{02}y_0y_2)+c_{01}x_1^2y_0y_1$.

\bf (iii) \rm Let $P=P_1\times P_2\in S_{\text{sing}}$ which is not an $A_1$-type and the fibre over $P_1$ is non-reduced. Then
from Theorem \ref{thm1.2}, the polynomial $f=a_{22}x_0^2y_2^2+x_0x_1(b_{11}y_1^2+b_{22}y_2^2+b_{02}y_0y_2+b_{12}y_1y_2)+x_1^2(c_{00}y_0^2+c_{11}y_1^2+c_{22}y_2^2+c_{01}y_0y_1+c_{02}y_0y_2+c_{12}y_1y_2)$, where $a_{22}\neq0$, $b_{11}\neq0$ and $c_{00}\neq0$. Let $\lambda(t)=\text{diag}(t^{-1},t)\times \text{diag}(t^{-1},t^0,t)$. Then $\lim\limits_{t\rightarrow 0} \lambda(t).f=a_{22}x_0^2y_2^2+x_0x_1(b_{11}y_1^2+b_{02}y_0y_2)+c_{00}x_1^2y_0^2$. 
Hence the result follows.
\end{proof}
\end{thm}
\section{Semi-stability and unstability of non-irreducible surfaces}
 In this section, we describe semi-stable and unstable (2,2)-type surfaces. We observe that non-irreducible (2,2)-type surfaces are never  stable.
Some arguments of proofs of this section are similar to the proofs of irreducible surfaces. So we skip them.
\begin{lem}
 Any surface of type (1,0), (0,1) and (1,1) of $\mathbb{P}^1\times\mathbb{P}^2$ is unstable by the natural linear action of the group $G=\text{SL}(2)\times \text{SL}(3)$ on corresponding linear systems. 
 \begin{proof}
  We know that $X^{ss}=\emptyset$, when $X=\mathbb{P}^n$ and the group acting on $X$ is SL$(n+1)$ (see Example 8.1 \cite{D}). Hence any surface $S$ either from $\vert(1,0)\vert$ or $\vert(0,1)\vert$ is unstable.
  
  Now let $S\in\vert(1,1)\vert$ and $S=Z(f)$. After some suitable linear change we can always assume $f=x_0y_2+x_1y_1$. Let $\lambda(t)=\text{diag}(t^{-1},t)\times \text{diag} (t^{-3},t,t^2)$. Then $\lim\limits_{t\rightarrow 0} \lambda(t).f=0$. Hence $S$ is unstable.
 \end{proof}
\end{lem}
\begin{lem}
 There is no stable surface in linear systems $\vert(2,0)\vert$ and $\vert(0,2)\vert$.
 \begin{proof}
  This result will follow from the fact that the space of quadric from $k[x_0,\cdots,x_n]_2$ has no stable points with respect to the action SL$(n+1)$ (see Example 10.1 \cite{D}).
 \end{proof}
\end{lem}
Let $S_1=Z(f_1)$ be an irreducible $(1,1)$-type surface of $\mathbb{P}^1\times\mathbb{P}^2$ and $i:S_1 \hookrightarrow \mathbb{P}^1\times\mathbb{P}^2$. Note that any irreducible $(1,1)$-type surface is always smooth. Let $p_1:S_1\rightarrow \mathbb{P}^1$ and  $p_2:S_1\rightarrow \mathbb{P}^2$ be natural projection maps. Also it can be checked easily that $S_1$ is isomorphic to $\mathbb{P}^2$ blown up at one point, $p_2$ is the blow-up map; and $S_1$ is isomorphic to the  Hirzebruch surface $\Sigma_1$, $p_1$ is the projection map.

Now let $Z(f_1), Z(f_2)\in \vert(1,1)\vert$ be irreducible surfaces and $f=f_1f_2$. Hence $Z(f)=Z(f_1)\cup Z(f_2)\in\vert(2,2)\vert$. Then the singularities of $S=Z(f)$ occur along the intersection $Z(f_1)\cap Z(f_2)$. So $S_{\text{sing}}\subset S$ is a divisor in $S$ which is linearly equivalent to $i^*\mathcal{O}_{\mathbb{P}^1\times\mathbb{P}^2}(1,1)$. Moreover,
$$S_{\text{sing}}\sim i^*\mathcal{O}_{\mathbb{P}^1\times\mathbb{P}^2}(1,1)\sim 2p_2^*\mathcal{O}_{\mathbb{P}^2}(1)-E\sim  2p_1^*\mathcal{O}_{\mathbb{P}^1}(1)+\mathcal{O}_{\Sigma_1}(1)$$
where $E$ is the exceptional curve of the blow up map $p_2$ and $\mathcal{O}_{\Sigma_1}(1)$ is the normalized section of $p_1$.

In the following lemma and corollary we prove the unstability conditions of $S$ which is a union of two (1,1)-type surfaces.
\begin{lem}\label{lem3.9}
 Let $Z(f_1), Z(f_2)\in \vert(1,1)\vert$ be irreducible surfaces and $f=f_1f_2$. Then $Z(f)\in\vert(2,2)\vert$.  $S=Z(f)$ is unstable if and only if $Z(f_1)$ and $Z(f_2)$ have common fibres over a same point $P_1\in\mathbb{P}^1$.
 \begin{proof}
Let $S=Z(f)=Z(f_1f_2)$ be an unstable surface and let
$f_1=x_0(\alpha_0y_0+\alpha_1y_1+\alpha_2y_2)+x_1(\beta_0y_0+\beta_1y_1+\beta_2y_2)$ and $f_2=x_0(\gamma_0y_0+\gamma_1y_1+\gamma_2y_2)+x_1(\delta_0y_0+\delta_1y_1+\delta_2y_2)$ be two irreducible polynomials.
Then

$f=x_0^2\Big(\alpha_0\gamma_0y_0^2+\alpha_1\gamma_1y_1^2+\alpha_2\gamma_2y_2^2+(\alpha_0\gamma_1+\alpha_1\gamma_0)y_0y_1+(\alpha_0\gamma_2+\alpha_2\gamma_0)y_0y_2+(\alpha_1\gamma_2+\alpha_2\gamma_1)y_1y_2\Big)+x_0x_1\Big((\alpha_0\delta_0+\beta_0\gamma_0)y_0^2+(\alpha_1\delta_1+\beta_1\gamma_1)y_1^2+(\alpha_2\delta_2+\beta_2\gamma_2)y_2^2+(\alpha_0\delta_1+\alpha_1\delta_0+\beta_0\gamma_1+\beta_1\gamma_0)y_0y_1+(\alpha_0\delta_2+\alpha_2\delta_0+\beta_0\gamma_2+\beta_2\gamma_0)y_0y_2+(\alpha_1\delta_2+\alpha_2\delta_1+\beta_1\gamma_2+\beta_2\gamma_1)y_1y_2\Big)+x_1^2\Big(\beta_0\delta_0y_0^2+\beta_1\delta_1y_1^2+\beta_2\delta_2y_2^2+(\beta_0\delta_1+\beta_1\delta_0)y_0y_1+(\beta_0\delta_2+\beta_2\delta_0)y_0y_2+(\beta_1\delta_2+\beta_2\delta_1)y_1y_2\Big).$

As $S$ is unstable, there exists a normalized one-parameter subgroup $\lambda(t)=\text{diag}(t^{r_0},t^{r_1})\times \text{diag} (t^{s_0},t^{s_1},t^{s_2})$ such that $\lim\limits_{t \to 0} \lambda(t).f=\lambda(t).\sum a_{\alpha\beta}x^{\alpha}y^{\beta}= 0$, i.e., $r\alpha+s\beta>0$ for all $a_{\alpha\beta}\neq 0$. 
 
Note that from the proof of Lemma \ref{lem2.17}, if $S$ is unstable, then the coefficients of $x_0^2y_0^2$, $x_0^2y_0y_1$, $x_0^2y_0y_2$, $x_0x_1y_0^2$ and $x_0x_1y_0y_1$ are zero. 
 The coefficients of $x_0^2y_0^2$, $\alpha_0\gamma_0=0$. Assume $\alpha_0=0$ but $\gamma_0\neq0$. The coefficient of  $x_0^2y_0y_1$, $\alpha_0\gamma_1+\alpha_1\gamma_0=0$. As $\gamma_0\neq0$ and $\alpha_0=0$, $\alpha_1=0$. The coefficient of  $x_0x_1y_0^2$, $\alpha_0\delta_0+\beta_0\gamma_0=0$. As $\gamma_0\neq0$ and $\alpha_0=0$, $\beta_0=0$. The coefficient of $x_0x_1y_0y_1$, $\alpha_0\delta_1+\alpha_1\delta_0+\beta_0\gamma_1+\beta_1\gamma_0=0$. As $\gamma_0\neq0$ and $\alpha_0=\alpha_1=\beta_0=0$, $\beta_1=0$. Hence we have $\alpha_0=\alpha_1=\beta_0=\beta_1=0$, which contradict the irreducibility of $f_1$. Hence our assumption was false. Similarly $\alpha_0 \neq 0$ but $\gamma_0=0$ contradicts the irreducibility of $f_2$. Therefore both $\alpha_0=\gamma_0=0$. Now
 
 $f=x_0^2\Big(\alpha_1\gamma_1y_1^2+\alpha_2\gamma_2y_2^2+(\alpha_1\gamma_2+\alpha_2\gamma_1)y_1y_2\Big)+x_0x_1\Big((\alpha_1\delta_1+\beta_1\gamma_1)y_1^2+(\alpha_2\delta_2+\beta_2\gamma_2)y_2^2+(\alpha_1\delta_0+\beta_0\gamma_1)y_0y_1+(\alpha_2\delta_0+\beta_0\gamma_2)y_0y_2+(\alpha_1\delta_2+\alpha_2\delta_1+\beta_1\gamma_2+\beta_2\gamma_1)y_1y_2\Big)+x_1^2\Big(\beta_0\delta_0y_0^2+\beta_1\delta_1y_1^2+\beta_2\delta_2y_2^2+(\beta_0\delta_1+\beta_1\delta_0)y_0y_1+(\beta_0\delta_2+\beta_2\delta_0)y_0y_2+(\beta_1\delta_2+\beta_2\delta_1)y_1y_2\Big).$
 
 As the coefficient of $x_0x_1y_0y_1$ is zero, $\alpha_1\delta_0+\beta_0\gamma_1=0$. As $S$ is unstable, either the coefficient of $x_0^2y_1^2$ or the coefficient of $x_1^2y_0^2$ is zero. Then we have either $\alpha_1\gamma_1=0$ or $\beta_0\delta_0=0$. Combining these three equations either $\alpha_1=0$ and $ \gamma_1=0$ or $\alpha_1=0$ and $\beta_0=0$. 
 
If  $\alpha_1=0$ and $ \gamma_1=0$, then
$$f=f_1f_2=\Big(\alpha_2x_0y_2+x_1(\beta_0y_0+\beta_1y_1+\beta_2y_2)\Big)\Big(\gamma_2x_0y_2+x_1(\delta_0y_0+\delta_1y_1+\delta_2y_2)\Big).$$
Hence $Z(f_1)$ and $Z(f_2)$ have same fibre over $[1,0]$, which is $[1,0]\times Z(y_2)$.
 
 If $\alpha_1=0$ and $\beta_0=0$, then
 
 $f=x_0^2\Big(\alpha_2\gamma_2y_2^2+\alpha_2\gamma_1y_1y_2\Big)+x_0x_1\Big(\beta_1\gamma_1y_1^2+(\alpha_2\delta_2+\beta_2\gamma_2)y_2^2+\alpha_2\delta_0y_0y_2+(\alpha_2\delta_1+\beta_1\gamma_2+\beta_2\gamma_1)y_1y_2\Big)+x_1^2\Big(\beta_1\delta_1y_1^2+\beta_2\delta_2y_2^2+\beta_1\delta_0y_0y_1+\beta_2\delta_0y_0y_2+(\beta_1\delta_2+\beta_2\delta_1)y_1y_2\Big).$
 
In this situation, $\alpha_2\neq0$ and $\beta_1\neq0$ as $f_1$ is irreducible. Also the unstability of $S$ implies that at least one of the coefficients of $x_0x_1y_1^2$ and $x_0x_1y_0y_2$ is zero (from Lemma \ref{lem2.17}). Hence either $\beta_1\gamma_1=0$ or $\alpha_2\delta_0=0$ i.e., $\gamma_1=0$ or $\delta_0=0$.
 
 Now if $\alpha_1=0$, $\beta_0=0$ and $\delta_0=0$, then
 $$f=f_1f_2=\Big(\alpha_2x_0y_2+x_1(\beta_1y_1+\beta_2y_2)\Big)\Big(x_0(\gamma_1y_1+\gamma_2y_2)+x_1(\delta_1y_1+\delta_2y_2)\Big).$$
 Note that $Z(f_1)$ and $Z(f_2)$ both have common section of $p_1$, which is $\mathbb{P}^1\times [1,0,0] \subset\mathbb{P}^1\times\mathbb{P}^2$. 
 Therefore $ \mathbb{P}^1\times [1,0,0]\subset S_{\text{sing}}$. But we observed that $S_{\text{sing}}\sim 2F+\sigma$ in $(1,1)$-type surface $Z(f_i)$, where $F$ is the numerical class of fibres and $\sigma$ is the normalized section of the projection map. 
 As the common section is a component of $S_{\text{sing}}$, common fibre is also a component of $S_{\text{sing}}$. So $Z(f_1)$ and $Z(f_2)$ have common fibres over some $P_1\in\mathbb{P}^1$. 
 
 Now we prove the converse part. Let $[1,0]\times [1,0,0]\in S_{\text{sing}}$. Then $f_1=x_0(\alpha_1y_1+\alpha_2y_2)+x_1(\beta_0y_0+\beta_1y_1+\beta_2y_2)$ and $f_2=x_0(\gamma_1y_1+\gamma_2y_2)+x_1(\delta_0y_0+\delta_1y_1+\delta_2y_2)$. 
 
 Let fibres over $[1,0]$ are same in $Z(f_1)$ and $Z(f_2)$. Then 
 $$f=\Big(\alpha_2x_0y_2+x_1(\beta_0y_0+\beta_1y_1+\beta_2y_2)\Big)\Big (\gamma_2x_0y_2+x_1(\delta_0y_0+\delta_1y_1+\delta_2y_2)\Big).$$
 Now let 
 $\lambda(t)=\text{diag}(t^{-3},t^3)\times \text{diag}(t^{-2},t^{-2},t^4)$. Then   $\lim\limits_{t\to 0}\lambda(t).f=0$.
 \end{proof}
\end{lem}
\begin{corl}\label{lem3.15}
 Let $Z(f_1), Z(f_2)\in \vert(1,1)\vert$ be irreducible polynomials and $f=f_1f_2$. Then $Z(f)\in\vert(2,2)\vert$. If  $S=Z(f)$ is not an unstable surface (i.e., each fibre over $[\alpha_0,\alpha_1]\in\mathbb{P}^1$ in $Z(f_1)$ and $ Z(f_2)$ are different lines in $\mathbb{P}^2$), then it is strictly semi-stable. Moreover, $f$ degenerates to $(\alpha_2x_0y_2+\beta_1x_1y_1)(\gamma_1x_0y_1+\delta_0x_1y_0)$, where $\alpha_2$, $\beta_1$, $\gamma_1$, and $\delta_0$ are non-zero.
 \begin{proof}
 Let $S=Z(f)=Z(f_1)\cup Z(f_2)$ such that $Z(f_1)$ and $Z(f_2)$ are $\mathbb{P}^2$ blow-up at $P_1$ and $P_2$ respectively. Then $\mathbb{P}^1\times P_1$ is the exceptional curve in $Z(f_1)$. The fibres of $Z(f_2)$ over $\mathbb{P}^1$ correspond to lines passing through the point $P_2$. Then there exists $Q_1\in\mathbb{P}^1$ such that the fibre over $Q_1$ is $Z(f_2)$ which is the line joining $P_1$ and $P_2$ i.e., $p_1^{-1}(Q_1)=Q_1\times L_{P_1P_2}$. Then $Q_1\times P_1\in Z(f_1)\cap Z(f_2)= S_{\text{sing}}$. So without loss of generality we may assume that $Q_1\times P_1=[1,0]\times [1,0,0]$ and $\mathbb{P}^1\times [1,0,0]$ is the exceptional curve of $Z(f_1)$. Then
 $$f=f_1f_2=\Big(\alpha_2x_0y_2+x_1(\beta_1y_1+\beta_2y_2)\Big)\Big(x_0(\gamma_1y_1+\gamma_2y_2)+x_1(\delta_0y_0+\delta_1y_1+\delta_2y_2)\Big).$$
 
Now let   $\lambda(t)=\text{diag}(t^{-1},t)\times \text{diag}(t^{-2},t^0,t^2)$. Then  $\lim\limits_{t\to 0}\lambda(t).f= (\alpha_2x_0y_2+\beta_1x_1y_1)(\gamma_1x_0y_1+\delta_0x_1y_0)$.
 Hence the result follows.
 \end{proof}
\end{corl}
\begin{prop}\label{prop3.15}
 Let $S=Z(f)$ be a non-irreducible (2,2)-type surface and $f=f_1f_2$, where $Z(f_1)\in \vert(1,0)\vert$, $Z(f_2)\in \vert(1,2)\vert$ and generic fibres of  $Z(f_2)$ over $\mathbb{P}^1$ are smooth. Then $S$ is semi-stable if and only if $C=Z(f_1)\cap Z(f_2)$ is smooth. Also $f$ degenerates to $x_0x_1(\alpha y_0y_2+\beta y_1^2)$, where $\alpha\neq 0$ and $\beta\neq0$.
 \begin{proof}
 Let $f=f_1f_2=(\gamma_0x_0+\gamma_1x_1)(x_0g_0(y_0,y_1,y_2)+x_1g_1(y_0,y_1,y_2))$, where $g_0(y_0,y_1,y_2)=\alpha_{00}y_0^2+\alpha_{11}y_1^2+\alpha_{22}y_2^2+\alpha_{01}y_0y_1+\alpha_{02}y_0y_2+\alpha_{12}y_1y_2$  and $g_1(y_0,y_1,y_2)=\beta_{00}y_0^2+\beta_{11}y_1^2+\beta_{22}y_2^2+\beta_{01}y_0y_1+\beta_{02}y_0y_2+\beta_{12}y_1y_2$ are two conics. 
 
 Assume that $f$ is unstable.  Then there exists a normalized one-parameter subgroup $\lambda(t)=\text{diag}(t^{r_0},t^{r_1})\times \text{diag} (t^{s_0},t^{s_1},t^{s_2})$ such that $\lim\limits_{t\to 0} \lambda(t).f=0$ where $r_0+r_1=0$, $r_0\leq r_1$;   $s_0+s_1+s_2=0$, $s_0\leq s_1\leq s_2$; and $r_0+s_0<0$. We know from the proof of Lemma \ref{lem2.17} that coefficients of $x_0^2y_0^2$, $x_0^2y_0y_1$ and $x_0^2y_0y_2$ are zero which imply $\alpha_{00}\gamma_0=0$, $\alpha_{01}\gamma_0=0$ and $\alpha_{02}\gamma_0=0$ respectively. 
 
 If $\gamma_0\neq0$, then $\alpha_{00}=\alpha_{01}=\alpha_{02}=0$. Note that coefficients of $x_0x_1y_0^2$ and $x_0x_1y_0y_1$ are also zero (from Lemma \ref{lem2.17}). Then $\beta_{00}=\beta_{01}=0$.
 Also unstability of $S$ implies that the coefficients of $x_0x_1y_0y_2$ and $x_0x_1y_1^2$; and  coefficients of $x_0x_1y_0y_2$ and $x_0^2y_1^2$ are not simultaneously non-zero. Therefore either $\beta_{02}=0$ or $\alpha_{11}=0$ and $\beta_{11}=0$.
 
 If $\beta_{02}=0$, then 
 $f=f_1f_2=(\gamma_0x_0+\gamma_1x_1)(x_0(\alpha_{11}y_1^2+\alpha_{22}y_2^2+\alpha_{12}y_1y_2)+x_1(\beta_{11}y_1^2+\beta_{22}y_2^2+\beta_{12}y_1y_2))$, where all fibres are non-smooth.
 
 If $\beta_{02}\neq0$, then $\alpha_{11}=\beta_{11}=0$. Therefore
 $f=f_1f_2=(\gamma_0x_0+\gamma_1x_1)(x_0(\alpha_{22}y_2^2+\alpha_{12}y_1y_2)+x_1(\beta_{22}y_2^2+\beta_{02}y_0y_2+\beta_{12}y_1y_2))$, where all fibres are non-smooth.
 
 If $\gamma_0=0$, then without loss of generality 
 $f=f_1f_2=x_1(x_0g_0(y_0,y_1,y_2)+x_1g_1(y_0,y_1,y_2))$. 
 Note that coefficients of $x_0x_1y_0^2$ and $x_0x_1y_0y_1$ are also zero (from Lemma \ref{lem2.17}). Hence $\alpha_{00}=0$ and $\alpha_{01}=0$. Then $g_0=\alpha_{11}y_1^2+\alpha_{22}y_2^2+\alpha_{02}y_0y_2+\alpha_{12}y_1y_2$.  If the coefficient of $x_0x_1y_1^2$, $\alpha_{11}\neq0$ and the coefficient of $x_0x_1y_0y_2$, $\alpha_{02}\neq0$, then $s_1>0$ and $s_0+s_2>0$ respectively which are not possible simultaneously from the condition of $\lambda$.
 So either $\alpha_{11}=0$ or $\alpha_{02}=0$. But for any of this cases  $g_0$ is non-irreducible. Hence $Z(f_1)\cap Z(f_2)=Z(g_0)$ is not smooth.
 
 Now we prove the converse part.
 Let $Z(f_1)\in \vert(1,0)\vert$, $Z(f_2)\in \vert(1,2)\vert$, $f=f_1f_2$; $C=Z(f_1)\cap Z(f_2)$. Then we prove if $C$ is non-smooth, then $Z(f)$ is unstable. After a suitable coordinate change, we may assume that $f_1=x_1$ and $f_2=x_0g_0(y_0,y_1,y_2)+x_1g_1(y_0,y_1,y_2)$. 
 From the assumption $C=Z(g_0)$  is a non-smooth curve. Then after some coordinate change
 $g_0(y_0,y_1,y_2)=y_2(\alpha_{02}y_0+\alpha_{12}y_1)$ or $g_0(y_0,y_1,y_2)=\alpha_{22}y_2^2$ when $g_0(y_0,y_1,y_2)$ is a reduced or non-reduced polynomial respectively. 
  
  Now let $\lambda(t)=\text{diag}(t^{-3},t^3)\times \text{diag}(t^{-2},t^{-2},t^4)$. Then $\lim\limits_{t\to0}\lambda(t).f=0$ and $f$ is unstable. Hence we complete the proof of the first part of the proposition.
 
 Let $f$ be a semi-stable polynomial.
 Now if $g_1$ is irreducible, then  after  some coordinate change we may assume $g_0=\alpha_{11}y_1^2+\alpha_{02}y_0y_2$, where $\alpha_{11}\neq 0$ and $\alpha_{02}\neq0$, and $f=f_1f_2=x_1(x_0(\alpha_{11}y_1^2+\alpha_{02}y_0y_2)+x_1g_1(y_0,y_1,y_2))$.
 
 Let  $\lambda(t)=\text{diag}(t^{-2},t^{2})\times \text{diag} (t^{-1},t^{0},t)$. Then $\lim\limits_{t\to 0} \lambda(t).f=x_1x_0(\alpha_{11}y_1^2+\alpha_{02}y_0y_2)$. Therefore $Z(f)$ is a strictly semi-stable surface. 
 \end{proof}
\end{prop}
In the following theorem we study all non-irreducible semi-stable $(2,2)$-type surfaces. We see that non-irreducible surfaces never be stable. Also all strictly semi-stable surfaces are either surfaces described in Corollary \ref{lem3.15} or surfaces of Proposition \ref{prop3.15}. 
 \begin{thm}\label{thm3.16}
Let $S=Z(f)$ be a non-irreducible (2,2)-type surface of $\mathbb{P}^1\times\mathbb{P}^2$. Then $S$ never be a stable surface.

$S$ is strictly semi-stable if and only if one of the followings holds:

(1) $f=f_1f_2$, where $Z(f_1)\in \vert(1,0)\vert$, $Z(f_2)\in \vert(1,2)\vert$, generic fibres of $Z(f_2)$ over $\mathbb{P}^1$ are smooth and $C=Z(f_1)\cap Z(f_2)$ is smooth.
 Also $f$ degenerates to $x_0x_1(\alpha y_0y_2+\beta y_1^2)$, $\alpha\neq 0$ and $\beta\neq 0$.
 
(2)   $f=f_1f_2$, $Z(f_1)\in \vert(1,1)\vert$, $Z(f_2)\in  \vert(1,1)\vert$ and  each fibre over $[\alpha_0,\alpha_1]\in\mathbb{P}^1$ in $Z(f_1)$ and $ Z(f_2)$ are different lines in $\mathbb{P}^2$. Moreover $f$ degenerates to $(\alpha_2x_0y_2+\beta_1x_1y_1)(\gamma_1x_0y_1+\delta_0x_1y_0)$.
 \begin{proof}
  Let $S=Z(f)$ be a non-irreducible $(2,2)$-type surface of $\mathbb{P}^1\times \mathbb{P}^2$. Then $f=f_1f_2$ where  $Z(f_1)\in \vert(1,0)\vert$ and $Z(f_2)\in \vert(1,2)\vert$, $Z(f_1)\in \vert(2,0)\vert$ and $Z(f_2)\in \vert(0,2)\vert$, $Z(f_1)\in \vert(1,1)\vert$ and $Z(f_2)\in \vert(1,1)\vert$ or $Z(f_1)\in \vert(2,1)\vert$ and $Z(f_2)\in \vert(0,1)\vert$.
  
  $\underline{\text{Case-1}:}$ $Z(f_1)\in \vert(1,0)\vert$ and $Z(f_2)\in \vert(1,2)\vert$.\\
  This case is studied in Proposition \ref{prop3.15} which are strictly semi-stable surfaces as described in (1) in this theorem.
 
$\underline{\text{Case-2}:}$
  $Z(f_1)\in \vert(2,0)\vert$ and $Z(f_2)\in \vert(0,2)\vert$\\
 After some coordinate change, $f_1$ is either $x_0x_1$ or $x_1^2$ when $f_1$ is  reduced and non-reduced respectively.  Up to some coordinate change  $f_2$ is  $a_{02}y_0y_2+a_{11}y_1^2$, $y_2(a_{02}y_0+a_{12}y_1)$ or $a_{22}y_2^2$ when $f_2$ is irreducible, reduced or non-reduced respectively.
 
 If $f_1$ is non-reduced, then $f$ is unstable with respect to $\lambda(t)=\text{diag}(t^{-3},t^3)\times \text{diag}(t^{-2},t^{-2},t^4)$. If $f_2$ is not irreducible, then $f$ is unstable and the one-parameter subgroup $\lambda(t)=\text{diag}(t^{-3},t^3)\times \text{diag}(t^{-2},t^{-2},t^4)$ is responsible for the unstability of $f$.

Now if $f=x_0x_1(a_{02}y_0y_2+a_{11}y_1^2)$ then from Case-1 $f$ is strictly semi-stable.

$\underline{\text{Case-3}:}$
  $Z(f_1)\in \vert(1,1)\vert$ and $Z(f_2)\in \vert(1,1)\vert$, where $f_1$ and $f_2$ are irreducible.\\
This case  follows from Lemma \ref{lem3.9} and Corollary \ref{lem3.15} which imply strictly semi-stable surfaces of (2) in this theorem.
 
$\underline{\text{Case-4}:}$
  $Z(f_1)\in \vert(2,1)\vert$ and $Z(f_2)\in \vert(0,1)\vert$\\
 After performing some linear changes, $f=y_2(y_0x_0x_1+y_1g_1(x_0,x_1)+y_2g_2(x_0,x_1))$. Let $\lambda(t)=\text{diag}(t^{-1},t)\times \text{diag}(t^{-3},t^{-1},t^4)$. Then $\lim\limits_{t\rightarrow 0}\lambda(t).f=0$. Hence $f$ is unstable.
 \end{proof}
\end{thm}
\section{Monomial reduction of strictly semi-stable and unstable (2,2)-type surfaces}
Here we generalize Mukai's idea (\cite{M}, Chapter 7.2) that  a normalized one-parameter subgroup of SL$(n+1)$ corresponds to two subsets of monomials which generate non-stable and unstable hypersurfaces respectively.

Let $M$ be the set of monomials which forms a basis  of the vector space $V=\mathbb{C}[x_0,x_1]_2\otimes\mathbb{C}[y_0,y_1,y_2]_2$. 
Given a normalized one-parameter subgroup\\
$\lambda(t)=\text{diag}(t^{r_0},t^{r_1})\times \text{diag}(t^{s_0},t^{s_1}, t^{s_2})$, one can define two subsets of $M$;

$ M^{+}(\lambda)=M^{+}((r_0,r_1),(s_0,s_1,s_2))=\{x^{\alpha}y^{\beta}=x_{0}^{\alpha_0}x_1^{\alpha_1}y_0^{\beta_0}y_1^{\beta_1}y_2^{\beta_2}\in M\vert r_0\alpha_0+r_1\alpha_1+s_0\beta_0+s_1\beta_1+s_2\beta_2>0\}$ and

$ M^{\oplus}(\lambda)=M^{\oplus}((r_0,r_1),(s_0,s_1,s_2))=\{x^{\alpha}y^{\beta}=x_{0}^{\alpha_0}x_1^{\alpha_1}y_0^{\beta_0}y_1^{\beta_1}y_2^{\beta_2}\in M\vert r_0\alpha_0+r_1\alpha_1+s_0\beta_0+s_1\beta_1+s_2\beta_2\geq0\}.$

If a (2,2)-type surface is unstable (resp. non-stable) with respect to $\lambda$, then the corresponding equation is a linear combination of monomials of $ M^{+}(\lambda)$ (resp. $M^{\oplus}(\lambda)$).

\begin{thm}
 Let $S=Z(h)$ be an  unstable (2,2)-type surface. Then there exists $f\in G.h$ such that $f$  is a combination  of $ M^{+}((-3,3),(-2,-2,4))$,  $ M^{+}((-4,4),(-10,5,5))$, $ M^{+}((-1,1),(-3,-1,4))$, or $ M^{+}((-2,2),(-5,-1,6))$ which are distinct subsets of $M$.
\begin{proof}
This follows from proofs of 
 Theorem \ref{thm1.1} and Theorem \ref{thm3.16}.
\end{proof}
\end{thm}
\begin{thm}\label{xrem3.20}
 Let $S=Z(h)$ be a strictly semi-stable (2,2)-type surface. Then there exists $f\in G.h$ such that $f$ is a combination  of $ M^{\oplus}((0,0),(-1,0,1))$,  $ M^{\oplus}((-1,1),(-2,0,2))$,  $ M^{\oplus}((-1,1),(-1,0,1))$ or $ M^{\oplus}((-2,2),(-1,0,1))$. Note that these are distinct subsets of $M$.
 \begin{proof}
The result clearly follows  from the proofs of  Theorem \ref{thm1.2} and Theorem \ref{thm3.16}(or Corollary \ref{lem3.15} and Proposition \ref{prop3.15}).
 \end{proof}
\end{thm}
\section{Compactification of the moduli space}

Let us denote by $\vert(2,2)\vert^{\text{s}}$ (resp.  $\vert(2,2)\vert^{\text{ss}}$)  the set of stable (resp. semi-stable) (2,2)-type surfaces by the $G=\text{SL}(2)\times \text{SL}(3)$ action. $\vert(2,2)\vert^{\text{s}}/G$ is called the moduli of (2,2)-type surface and $\vert(2,2)\vert^{\text{ss}}//G$ is it's natural compactification which is a projective variety.  The categorical quotient
 $$\vert(2,2)\vert^{\text{ss}}//G\cong \text{Proj}(R^G)$$ 
 where $R=\oplus_{n\geq0} H^0(\mathbb{P}^{17},\mathcal{O}_{\mathbb{P}^{17}}(n))$ and the geometric quotient $\vert(2,2)\vert^{\text{s}}/G$ is an open subset of $\text{Proj}(R^G)$. As the orbit of the strictly semi-stable surfaces are in the boundary of the compactification, maximal elements of sets $ M^{\oplus}((r_0,r_1),(s_0,s_1,s_2))$ (see Theorem \ref{xrem3.20}) correspond to boundary elements. If $S\in \vert(2,2)\vert$, then $[S]$ is the  orbit in $\vert(2,2)\vert//G$. In the following theorem, we calculate the dimension of the moduli and describe the behaviour of boundary stratum of the moduli.

 \subsection{Proof of Theorem \ref{thm1.3}}
The dimension of the moduli space is given by
$$\text{dim}\Big(\vert(2,2)\vert^{\text{s}}/G\Big)=\text{dim}(\vert(2,2)\vert)-\text{dim}(G)=17-11=6.$$
Let $Z(f')$ be a strictly semi-stable (2,2)-type surface of $\mathbb{P}^1\times\mathbb{P}^2$. Then there exists a one-parameter subgroup $\lambda'$, such that $\lim\limits_{t\rightarrow 0}\lambda'(t).f'=h'\neq 0$. As $\lambda'$ is  diagonalizable for some $g\in G$, $\lambda'(t)=g^{-1}\lambda(t)g$, where $\lambda$ is a normalized one-parameter subgroup. So after a coordinate change with respect to $g$, we have  $\lim\limits_{t\rightarrow 0}\lambda(t).f=h\neq 0$, where
 $g.f'=f$ and $g.h'=h$.
We know that if two orbit closure of two distinct points meet then their images in the quotient space are same  (Corollary 6.1 in \cite{D}).  
Then from Theorem \ref{thm3.11} and Theorem \ref{thm3.16}, following are the sets of strictly semi-stable surfaces:\\
 $\Gamma_1=\Big\{Z(x_1x_0(b_{02} y_0y_2+b_{11} y_1^2))\mid$ $b_{02}\neq0$ and $b_{11}\neq0\Big\}$,\\
 $\Gamma_2=\Big\{Z(a_{12}x_0^2y_1y_2+x_0x_1(b_{11}y_1^2+b_{02}y_0y_2)+c_{01}x_1^2y_0y_1)\mid$ either $b_{02}\neq0$, $a_{12}\neq0$, and $c_{01}\neq0$ or $b_{02}\neq0$ and $b_{11}\neq0\Big\}$,\\
 $\Gamma_3=\Big\{Z(a_{22}x_0^2y_2^2+x_0x_1(b_{11}y_1^2+b_{02}y_0y_2)+c_{00}x_1^2y_0^2)\mid$ either $b_{11}\neq0$, $a_{22}\neq 0$ and $c_{00}\neq0$ or  $b_{11}\neq0$ and  $b_{02}\neq0\Big\}$,\\
 $\Gamma_4=\Big\{Z(x_0^2(a_{11}y_1^2+a_{02}y_0y_2)+x_0x_1(b_{11}y_1^2+b_{02}y_0y_2)+x_1^2(c_{11}y_1^2+c_{02}y_0y_2))\mid$  $a_{11}\neq0$, $c_{11}\neq0$ and $b_{02}\neq0$; $a_{02}\neq0$, $c_{02}\neq0$ and $b_{11}\neq0$; or $b_{11}\neq0$ and $b_{02}\neq0\Big\}$. 
 
 We denote the images of $\Gamma_i$ by $\overline{\Gamma_i}$ in  $\vert(2,2)\vert^{\text{ss}}//G$ for $i=1,2,3,4$.
 
\bf We claim that \rm $\overline{\Gamma_1}$ is a point and  $\overline{\Gamma_i}$, for $i=2,3,4$ are rational curves in  $\vert(2,2)\vert^{\text{ss}}//G$. Furthermore,  we also claim $\overline{\Gamma_1}\in\overline{\Gamma_i}$, for $i=2,3,4$.

 Note that $g.(x_1x_0(b_{02} y_0y_2+b_{11} y_1^2))=b_{11}^{1/3}b_{02}^{2/3}x_1x_0( y_0y_2+ y_1^2)$,
 
 where $g=\text{id}\times\text{diag}(b_{11}^{1/6}b_{02}^{-1/6},b_{11}^{-1/3}b_{02}^{1/3},b_{11}^{1/6}b_{02}^{-1/6})$.
 
 So $[Z(x_1x_0(b_{02} y_0y_2+b_{11} y_1^2))]=[Z(x_1x_0( y_0y_2+ y_1^2))]$ in $\overline{\Gamma_1}$. Hence $\overline{\Gamma_1}=[Z(x_1x_0( y_0y_2+ y_1^2))]$ is a point.
 
 Now consider the set of strictly semi-stable surfaces $\Gamma_2=\Big\{Z(a_{12}x_0^2y_1y_2+x_0x_1(b_{11}y_1^2+b_{02}y_0y_2)+c_{01}x_1^2y_0y_1)\mid$ either $b_{02}\neq0$, $a_{12}\neq0$, and $c_{01}\neq0$ or $b_{02}\neq0$ and $b_{11}\neq0\Big\}$.
 
 Note that $x_0x_1(b_{11}y_1^2+b_{02}y_0y_2)+c_{01}x_1^2y_0y_1=g.(a_{12}x_0^2y_1y_2+x_0x_1(b_{11}y_1^2+b_{02}y_0y_2))$ for some $g\in G$. So without loss
  of generality we may assume that $a_{12}\neq0$. Moreover $b_{02}\neq 0$, otherwise $f=a_{12}x_0^2y_1y_2+x_0x_1(b_{11}y_1^2+b_{02}y_0y_2)+c_{01}x_1^2y_0y_1$ becomes unstable. Therefore 
  $$g.(a_{12}x_0^2y_1y_2+x_0x_1(b_{11}y_1^2+b_{02}y_0y_2)+c_{01}x_1^2y_0y_1)=(\delta x_0^2y_1y_2+x_0x_1(b_{11}'y_1^2+\delta y_0y_2)+c_{01}'x_1^2y_0y_1),$$
  where $g=\text{id}\times\text{diag}(a_{12}b_{02}^{-1},1,a_{12}^{-1}b_{02})$. Without loss of generality, we may assume $\delta =1$. Now
  $$\overline{\Gamma_2}=\Big\{[Z(x_0^2y_1y_2+x_0x_1(b_{11}y_1^2+y_0y_2)+c_{01}x_1^2y_0y_1)]\mid  b_{11}\neq0 \text{ or } c_{01}\neq0\Big\}.$$
 Let us define a surjective map $\phi:\mathbb{A}^2\backslash \{(0,0)\}\rightarrow \overline{\Gamma_2}$ such that $\phi(u,v)=[Z(x_0^2y_1y_2+x_0x_1(uy_1^2+y_0y_2)+vx_1^2y_0y_1)]$. We claim that $\phi(u,v)=\phi(ru,rv)$, for any $r\in \mathbb{C}^*$.
 
 $g.(x_0^2y_1y_2+x_0x_1(ruy_1^2+y_0y_2)+rvx_1^2y_0y_1)=r^{1/3}(x_0^2y_1y_2+x_0x_1(uy_1^2+y_0y_2)+vx_1^2y_0y_1)$, where $g=\text{diag}(r^{1/4},r^{-1/4})\times \text{diag}(r^{1/6},r^{-1/3},r^{1/6})$. Hence we proved our claim and defined a morphism $\overline{\phi}:\mathbb{P}^1\rightarrow \overline{\Gamma_2}$. Therefore $\overline{\Gamma_2}$ is a rational curve. Moreover, $\overline{\phi}([1,0])=\overline{\Gamma_1}$ and $\overline{\phi}([0,1])=[Z(x_0^2y_1y_2+x_0x_1y_0y_2+x_1^2y_0y_1)]$.
 
 
 Let us consider the set  $\Gamma_3=\Big\{Z(a_{22}x_0^2y_2^2+x_0x_1(b_{11}y_1^2+b_{02}y_0y_2)+c_{00}x_1^2y_0^2)\mid$ either $b_{11}\neq0$, $a_{22}\neq 0$ and $c_{00}\neq0$; or  $b_{11}\neq0$ and  $b_{02}\neq0\Big\}$ of strictly semi-stable surfaces.
 
 Note that $[Z(x_0x_1(b_{11}y_1^2+b_{02}y_0y_2)+c_{00}x_1^2y_0^2)]=[Z(a_{22}x_0^2y_2^2+x_0x_1(b_{11}y_1^2+b_{02}y_0y_2))]$ and $b_{11}\neq 0$ for any elements of $\Gamma_3$. Similarly, we can write 
 $$\overline{\Gamma_3}=\Big\{[Z(x_0^2y_2^2+x_0x_1(y_1^2+b_{02}y_0y_2)+c_{00}x_1^2y_0^2)]\mid  b_{02}\neq0 \text{ or } c_{00}\neq0\Big\}.$$
 Now we define a morphism $\phi:\mathbb{A}^1\rightarrow \overline{\Gamma_3}\backslash \{[Z(x_0^2y_2^2+x_0x_1y_1^2+x_1^2y_0^2)]\}$, such that $\phi(u)=[Z(x_0^2y_2^2+x_0x_1(y_1^2+y_0y_2)+ux_1^2y_0^2)]$. We claim that this $\phi$ is surjective. This follows immediately from 
 $g.(x_0^2y_2^2+x_0x_1(y_1^2+b_{02}y_0y_2)+c_{00}x_1^2y_0^2)=b_{02}^{2/3}(x_0^2y_2^2+x_0x_1(y_1^2+y_0y_2)+c_{00}b_{02}^{-2}x_1^2y_0^2)$, where $g=\text{id} \times \text{diag}(b_{02}^{-2/3}, b_{02}^{1/3},b_{02}^{1/3})$. Therefore we can extend the morphism $\overline{\phi}:\mathbb{P}^1\rightarrow \overline{\Gamma_3}$ which is also surjective. Moreover, $\overline{\phi}([1,0])=\overline{\Gamma_1} $ and $\overline{\phi}([0,1])=[Z(x_0^2y_2^2+x_0x_1y_1^2+x_1^2y_0^2)]$. Hence $\overline{\Gamma_3}$ is  a rational curve.
 
Finally, consider the set $\Gamma_4=\Big\{Z(x_0^2(a_{11}y_1^2+a_{02}y_0y_2)+x_0x_1(b_{11}y_1^2+b_{02}y_0y_2)+x_1^2(c_{11}y_1^2+c_{02}y_0y_2))\mid$  $a_{11}\neq0$, $c_{11}\neq0$ and $b_{02}\neq0$;  $a_{02}\neq0$, $c_{02}\neq0$ and $b_{11}\neq0$; or $b_{11}\neq0$ and $b_{02}\neq0\Big\}$.
 
 Then using a similar calculation, we can write the image of $\Gamma_4$ in $\vert(2,2)\vert^{\text{ss}}//G$ as
 $$\overline{\Gamma_4}=\Big\{[Z(x_0^2y_0y_2+x_0x_1(y_1^2+b_{02}y_0y_2)+x_1^2c_{02}y_0y_2)]\mid b_{02}\neq0;\text{ or }  c_{02}\neq0\Big\}$$
 and define a morphism $\phi:\mathbb{A}^1\rightarrow \overline{\Gamma_4}\backslash [Z(x_0^2y_0y_2+x_0x_1y_1^2+x_1^2c_{02}y_0y_2)]$, which can be extended to $\overline{\phi}:\mathbb{P}^1\rightarrow\overline{\Gamma_4}$.
 Therefore, $\overline{\Gamma_4}$ is a  rational curve and $\overline{\Gamma_1}\in \overline{\Gamma_4}$. Hence we prove our claims and the result follows.
 

\end{document}